\documentclass[ba]{imsart}
\pubyear{2025}
\volume{TBA}
\issue{TBA}
\firstpage{1}
\lastpage{1}

\usepackage{placeins}
\usepackage{pgfplots}
\usepackage{tikz}
\usepackage{pgfplotstable}
\pgfplotsset{width=10cm,compat=1.9}

\usepackage{amsfonts,amsmath,amssymb,amsthm}
\usepackage{natbib}
\usepackage{bm}
\usepackage{enumitem}
\usepackage{dsfont}
\usepackage{graphicx}
\usepackage{mathtools}
\usepackage[titlenumbered,ruled]{algorithm2e}
\usepackage{subcaption}
\usepackage{xcolor}
\usepackage{float}

\definecolor{coolgrey}{rgb}{0.55, 0.57, 0.67}
\definecolor{slategray}{rgb}{0.44, 0.5, 0.56}

\newtheorem{theorem}{Theorem}[section]
\newtheorem{proposition}{Proposition}[section]
\newtheorem{remark}{Remark}[section]
\newtheorem{assumption}{Assumption}
\newtheorem{lemma}[theorem]{Lemma}
\newtheorem{refproof}{Proof}

\newcommand{\dd}[0]{\mathrm{d}}
\newcommand{\bbE}[0]{\mathbb{E}}
\newcommand{\bbR}[0]{\mathbb{R}}
\newcommand{\bbV}[0]{\mathbb{V}}
\newcommand{\bbC}[0]{\mathbb{C}}
\newcommand{\calL}[0]{\mathcal{L}}
\newcommand{\calN}[0]{\mathcal{N}}
\newcommand{\wtX}[0]{\widetilde{X}}

\DeclarePairedDelimiter\norm{||}{||}
\DeclarePairedDelimiter\normLtwo{||}{||_{\calL^2}}
\DeclarePairedDelimiter\biggnormLtwo{\bigg|\bigg|}{\bigg|\bigg|_{\calL^2}}
\DeclarePairedDelimiter\bignormLtwo{\big|\big|}{\big|\big|_{\calL^2}}
\DeclarePairedDelimiter\bignorm{\bigg\|}{\bigg\|}
\DeclarePairedDelimiter\abs{|}{|}

\usepackage[colorlinks,citecolor=blue,urlcolor=blue,filecolor=blue,backref=page]{hyperref}

\newcommand{\yvann}[1]{\textcolor{blue}{YLF: #1}}

\startlocaldefs
\endlocaldefs

\begin{document}

    \begin{frontmatter}
        \title{Modelling pathwise uncertainty of Stochastic Differential Equations samplers via Probabilistic Numerics}
        \runtitle{Probabilistic Numerics for SDEs}

        \begin{aug}
            \author{\fnms{Yvann} \snm{Le Fay}\thanksref{ENSAE,Aalto}\ead[label=e1]{yvann.lefay@ensae.fr}
            },
            \author{\fnms{Simo} \snm{S\"{a}rkk\"{a}}\thanksref{Aalto}\ead[label=e2]{simo.sarkka@aalto.fi}}
            \and
            \author{\fnms{Adrien} \snm{Corenflos}\thanksref{Aalto,Warwick}\ead[label=e3]{adrien.corenflos@warwick.ac.uk}}
            \address[ENSAE]{ENSAE Paris, CREST, Institut Polytechnique de Paris\printead[presep={,\ }]{e1}}
            \address[Aalto]{Department of Electrical Engineering and Automation, Aalto University\printead[presep={,\ }]{e2}}
            \address[Warwick]{Department of Statistics, University of Warwick\printead[presep={,\ }]{e3}}
            \runauthor{Y. Le Fay, S. S\"arkk\"a, and A. Corenflos}

        \end{aug}
        \begin{abstract}
            Probabilistic ordinary differential equation (ODE) solvers have been introduced over the past decade as uncertainty-aware numerical integrators.
            They typically proceed by assuming a functional prior to the ODE solution, which is then queried on a grid to form a posterior distribution over the ODE solution.
            As the queries span the integration interval, the approximate posterior solution then converges to the true deterministic one.
            Gaussian ODE filters, in particular, have enjoyed a lot of attention due to their computational efficiency, the simplicity of their implementation, and their provable fast convergence rates.
            In this article, we extend the methodology to stochastic differential equations (SDEs) and propose a probabilistic simulator for SDEs.
            Our approach involves transforming the SDE into a sequence of random ODEs using piecewise differentiable approximations of the Brownian motion.
            We then apply probabilistic ODE solvers to the individual ODEs, resulting in a pathwise probabilistic solution to the SDE\@.
            We establish worst-case strong $1.5$ local and $1.0$ global convergence orders for a specific instance of our method.
            We further show how we can marginalise the Brownian approximations, by incorporating its coefficients as part of the prior ODE model, allowing for computing exact transition densities under our model.
            Finally, we numerically validate the theoretical findings, showcasing reasonable weak convergence properties in the marginalised version.
        \end{abstract}


    \end{frontmatter}

    \section{Introduction}\label{section:introduction}
    Consider the following It\^{o} stochastic differential equation (SDE) with an initial value condition and an additive noise:
    \begin{equation}
        \begin{split}
            \label{eq:additive_sde}
            \dd X_{t} &= \mu(X_{t}, t)\dd t + \sigma(t)\dd B_{t}, \indent t \in [0, T]\\
            X_{0} &= Z,
        \end{split}
    \end{equation}
    where $T > 0$ is a finite deterministic horizon and $\{ B_{t}~:~t \in [0,T]\}$ is an $m$-dimensional Brownian motion defined on some probability space, and $Z$ be a square-integrable random variable independent of the canonical filtration generated by $B$.
    Additionally, we assume the drift vector function $\mu:\bbR^d\times [0,T] \to \bbR^d$ and the diffusion matrix function $\sigma:[0,T] \to \bbR^{d \times m}$ are well-behaved functions, ensuring the existence of a unique strong solution $X$ to~\eqref{eq:additive_sde}~\citep[see, e.g.,][Ch. 5.3]{OK2003}.

    Given a discretisation $0=t_{0} < t_{1} < t_{2} < \cdots < t_{K} = T$, sampling from the finite-dimensional distribution $(X_{t_{0}}, \ldots X_{t_{K}})$ can then be done by sampling iteratively from the conditional densities $p_{t_{k}, t_{k+1}}(X_{t_{k+1}} | X_{t_{k}})$.
    However, closed-form or otherwise tractable expressions for these densities are rarely known, and we need instead to resort to numerical approximations.
    The simplest such approach is given by the Euler--Maruyama (EM) scheme~\citep[see, e.g.,][Ch. 10.2]{kloeden1992numericalSDE}, in which we informally speaking replace infinitesimal differences with finite differences, leveraging the fact that Brownian increments can be simulated exactly.
    Thus, the solution is approximated by
    \begin{align}
        \label{eq:euler_maruyama}
        X_{t_{k+1}} \approx X_{t_{k}} + \mu(X_{t_{k}}, t_{k}) (t_{k+1} - t_{k}) +  \sigma(t_{k}) (B_{t_{k+1}} - B_{t_{k}}).
    \end{align}
    The EM scheme~\eqref{eq:euler_maruyama} has both weak and strong convergence rates of order $1$, which means that a discrete-time approximation $X^{\delta}$, where $\delta$ is the maximum of the differences $t_{k+1} - t_{k}$, converges to the true solution $X$ as $\delta$ goes to $0$. Formally, it converges weakly with order $\beta > 0$, if, for any polynomial function $g$, $\max_{k\in [0, K]}\norm{\bbE(g(X_{t_{k}}))-\bbE(g(X^{\delta}_{t_{k}}))} = O(\delta^\beta)$, and it converges strongly with order $\gamma > 0$, if $\max_{k\in[0, K]} \bbE(\norm{X_{t_{k}} - X^{\delta}_{t_{k}}}) = O(\delta^\gamma)$~\citep[see, e.g.,][Ch. 9.6, Ch. 9.7]{kloeden1992numericalSDE}.

    Higher-order schemes~\citep[for a review, see, e.g.,][Ch. 10, Ch. 11]{kloeden1992numericalSDE} may sometimes be necessary, either for computational reasons (allowing for better precision at a similar cost) or for physical/statistical reasons.
    For example, when simulating stochastic Hamiltonian dynamics~\citep[see, e.g.,][]{Burrage2014}, energy preservation of the system is of paramount importance.
    The need for higher-order schemes is sometimes more subtle, as in \citet{Ditlevsen_2018} who demonstrate that, for instance, the EM scheme fails to capture the covariance of the transition density arising from hypoelliptic diffusions, preventing parameter estimation for such SDEs if a higher-order scheme is not used.

    The discretisation errors of the numerical approximation schemes add up a layer of additional numerical uncertainty on top of the desirable structural randomness, that is, the law of the process $X$ one would be interested in recovering.
    In this article, we aim to provide a model to quantify this numerical uncertainty by extending the probabilistic numerics~\citep[PN, ][]{hennig_osborne_kersting_2022} framework for ODEs to SDEs.
    At the core, our method relies on modelling these uncertainties in terms of a corrected posterior distribution over the solution.
    This perspective has recently been introduced for ordinary differential equations~\citep[ODEs, see, e.g.,][]{Chkrebtii2016BayesianODE,HennigKerstingTronarpSarkka2019,HennigTronarpSarkka2021} as part of the probabilistic numerics framework which we briefly review in Section~\ref{section:pn_for_odes}. Detailed contributions of the present article are given at the end of this section.

    \label{subsection:probnum-sde}
    This problem of quantifying uncertainty arising from discretisation in SDEs was raised by~\citet{Lysy2016Comment} in a commentary to the article of~\citet{Chkrebtii2016BayesianODE}, treating ODEs. \citet{Lysy2016Comment} highlights that the non-differentiability of the Brownian path hinders the use of the Bayesian methodology given in~\citet{Chkrebtii2016BayesianODE} which relies on jointly modelling the ODE solution and its first derivative by a Gaussian process.
    \citet{Lysy2016Comment} suggested two alternative approaches. The first suggestion is to approximate the Brownian motion by its pre-generated increments in an Euler--Maruyama fashion. The obtained algorithm is then an inference algorithm on the increments of the process. However, it leads to asymptotically biased drift and diffusion terms as the discretisation mesh becomes thinner, even for constant drift and diffusion functions.
    The second alternative is to replace the Brownian motion with a well-chosen differentiable approximating process. Specifically, \citet{Lysy2016Comment} proposes to use the integrated Ornstein-Uhlenbeck process (IOUP) $G_t$, a stationary Gaussian process with $\textup{cov}(\dot{G}_s, \dot{G}_t) = e^{-\theta |t - s|}$, after which~\eqref{eq:additive_sde} becomes
    \begin{equation*}
        \dd X_t = \mu(X_t, t)\dd t + \sigma(t)\dd{G}_t.
    \end{equation*}
    This approach resembles that of Wong--Zakai approximations~\citep{wongzakai65,Twardowska1996WongZakaiAF}, but also presents the same issue: for a choice of $\theta > 0$, the solution will exhibit a non-vanishing bias and, as $\theta \to 0$, the resulting ODE will become stiffer and, hence, harder to solve.
    Moreover, globally uniform Wong--Zakai approximations are, in general, computationally expensive for low convergence rates~\citep{GyongyWongZakaiConvergenceRates}.

    As stated in \citet{Lysy2016Comment}, a differentiable approximation of the Brownian motion, given that the associated random ODE recovers~\eqref{eq:additive_sde} in a suitable limit, enables us to apply the probabilistic numerics framework to SDEs. In this article, we use a piecewise differentiable approximation of the Brownian motion for this purpose.

    \paragraph{Outline and Contributions.}
    \label{section:outline}
    \begin{enumerate}
        \item \textbf{Goal: } Our manuscript aims to quantify pathwise uncertainty in SDEs with additive noise.
        To this end, we generalise the Bayesian probabilistic numerics framework for ODEs to handle SDEs with additive noise.
        \item \textbf{Method Overview: } Our method relies on (i) assuming a prior distribution over the possible paths of the SDE, like~\citet{Lysy2016Comment}, and (ii) deriving the posterior distribution thereof associated with a likelihood model corresponding to the path satisfying the SDE.
        \begin{enumerate}[label=(\roman*)]
            \item \textbf{Piecewise Differentiable Brownian Motion Approximation: } Instead of using a global approximation (e.g., the Wong-Zakai approximation), our approach uses piecewise differentiable approximations of the Brownian motion.
            \item \textbf{Bayesian Filtering and Smoothing: } Our inference process is framed within the Bayesian filtering and smoothing framework~\citep{sarkka2023bayesian}, which is commonly used in probabilistic numerics for ODEs~\citep[see, e.g.,][]{HennigKerstingTronarpSarkka2019}.
        \end{enumerate}
    \end{enumerate}

    In summary, our approach builds on Bayesian numerics for ODEs, extending it to SDEs with additive noise, using piecewise approximations of Brownian motion and probabilistic filtering for uncertainty quantification.

    The remainder of the article is organised as follows.
    In Section~\ref{section:background}, we review the piecewise approximation approach to SDE solving, in particular the method of~\citet{Foster2020OptimalPolyBrownian}, as well as probabilistic ODE solver methods with a focus on the Gaussian filters as given in~\citet{HennigKerstingTronarpSarkka2019}.
    In Section~\ref{section:method}, we show how probabilistic ODE solvers can be applied to random ODE approximations of SDEs. In the specific case of Gaussian filtering applied to the method of~\citet{Foster2020OptimalPolyBrownian}, we derive theoretical guarantees in the additive noise setting for the strong convergence of our method.
    Extensions of our proposed method are discussed in Sections~\ref{section:extension} and~\ref{section:extension2}.
    Most importantly, we discuss the explicit marginalisation of some auxiliary Gaussian variables introduced as the coefficients of the Brownian approximation.
    In Section~\ref{section:calibration}, we show how posterior uncertainty calibration can be achieved \emph{after the fact} via maximum likelihood estimation for the prior diffusion coefficients
    Finally, in Section~\ref{section:experiments}, we provide numerical evidence showcasing the statistical performances of our method and its extensions.

    \section{Background}
    \label{section:background}
    We now turn to review the two main ingredients of our method, the Bayesian framework to quantify numerical uncertainties for ODEs and the piecewise differentiable approximation of the Brownian motion we use to transform the SDE into a sequence of random ODEs. In Section~\ref{section:method}, we will then show how to principally combine Sections~\ref{section:gaussian_ode_filters} and~\ref{section:brownian_approximation} to perform pathwise uncertainty modelling in SDEs.

    \subsection{Bayesian Uncertainty Quantification for ODEs}
    \label{section:pn_for_odes}

    Consider the following initial value problem (IVP),
    \begin{equation}
        \begin{split}
            \label{eq:IVP}
            y'(t) &= f(y(t), t),  \quad \forall t\in [0, T] \\
            y(0) &= y_{0},
        \end{split}
    \end{equation}
    where $f:\bbR^{d}\times [0, T]\to \bbR^{d}$ is the vector field. %
    Solutions to~\eqref{eq:IVP} are rarely expressible in closed form, and one needs to resort to some numerical schemes.
    Traditional numerical integrators for ODEs such as Euler schemes or Runge-Kutta methods are solely point estimators of the trajectory.
    They typically come with a numerical error that takes the form $\norm{y^{h}-y}_{\infty}=O(h^{\alpha})$ where $y^h$ is the approximated solution with step size $h$.
    Different approaches to performing numerical analysis of the error exist, and recently proposed methods are given, for example, as randomised integration grids~\citep{Lie2022}, black-box error inference via interpolation of pre-existing numerical methods~\citep{teymur2021black}, or designing numerical methods that intrinsically carry uncertainty via Bayesian procedures, the origin of which can be traced back to at least~\citet{Skilling1992}.

    A recent outlook~\citep{HennigKerstingTronarpSarkka2019} to the Bayesian procedure consists in casting the IVP into a Bayesian inference problem in which $y$ is replaced by a random process $Y$, and is constrained using the measurement process:
    \begin{align}
        \label{eq:misalignement_measure}
        \mathcal{Z}(t) = Y^{(1)}(t) - f(Y^{(0)}(t), t).
    \end{align}
    Indeed, in the limiting case when $\mathcal{Z} \equiv 0$, the solution contracts to the true solution $y$ of~\eqref{eq:IVP},~\citep[see, e.g.,][Ch. 6., Sec. 38.1.1]{hennig_osborne_kersting_2022}.
    However, the non-linear dependence of $\mathcal{Z}$ on $Y$ hinders the exact computation of the posterior.
    Hence, it becomes necessary to resort to some approximation schemes for $\mathcal{Z}$ which are typically sample-based~\citep{Chkrebtii2016BayesianODE,HennigKerstingTronarpSarkka2019} or Gaussian-based~\citep{HennigKerstingTronarpSarkka2019,HennigTronarpSarkka2021}.
    In~\citet{Chkrebtii2016BayesianODE}, the prior-to-posterior update is obtained by sampling from the posterior predictive distribution and comparing it to the ODE at hand.
    This approach leads to expensive routines, typically cubic with respect to the number of time steps, but presents the benefit of performing ``exact'' inference in the limit of an infinite number of samples.
    On the other hand,~\citet{HennigTronarpSarkka2021} assume a Gauss-Markov process prior for $Y$, and approximate the posterior distribution through the use of Gaussian (Kalman) filters~\citep[see, e.g.,][]{sarkka2023bayesian}.
    In this case, the inference can be performed with a linear complexity in the size of the mesh, but results in biased (Gaussian) posterior representations.
    The specific merits of the two approaches are thoroughly discussed in the textbook~\citet[Ch. 6., Sec. 38.3.1, 38.4]{hennig_osborne_kersting_2022}, to which we refer the reader.
    Throughout this article, we consider the latter method, which we review in the next section, and refer to \citet{henning2015probs,Cockayne2019review}, and references within, for a historical review of the field.

    \subsection{The Gaussian ODE Filters}
    \label{section:gaussian_ode_filters}
    \citet{HennigKerstingTronarpSarkka2019,HennigTronarpSarkka2021} propose to use non-linear Gaussian filters to solve the IVP~\eqref{eq:IVP}.
    We quickly review their approach as well as existing convergence guarantees.
    Without loss of generality, we consider the unidimensional case $d=1$.
    As mentioned in Section~\ref{section:pn_for_odes}, conditioning $Y$ on $\mathcal{Z}=0$ on the whole interval $[0, T]$ contracts the posterior to the true solution of~\eqref{eq:IVP} but is intractable in practice.
    Consequently, we need to restrict ourselves to a finite conditioning grid.
    Let $0=t_{0}<t_{1}<\ldots<t_K=T$ be a regular discrete grid with step size $h$.
    The prior on $Y$ is set to be the solution of the following SDE,
    \begin{equation*}
        \begin{split}
            \dd Y(t) &= F Y(t)\dd t+ L \dd B_t, \quad \forall t\in [0, T]\\
            Y(0) &= \calN(y_{0}, P(0)),
        \end{split}
    \end{equation*}
    where the drift and diffusion matrices are given by
    \begin{equation}
        \label{eq:transition_matrix_ekf}
        F=\begin{pmatrix}
              0      & 1      & 0      & \ldots & 0      \\
              \vdots &        & \ddots &        & \vdots \\
              \vdots &        &        & \ddots & \vdots \\
              0      & \ldots & \ldots & 0      & 1      \\
              c_{0}  & \ldots & \ldots & \ldots & c_q
        \end{pmatrix},
        \indent L = \begin{pmatrix}
                        0 \\ \vdots \\\eta
        \end{pmatrix},
    \end{equation}
    for some fixed parameters $(c_{0}, \ldots, c_q, \eta)\in \bbR^{q+2}$.
    The solution of this SDE is a Gauss-Markov process satisfying, for any $t, h\geq 0$, $Y(t+h) \sim \calN(A(h)Y(t), Q(h))$, with $A(h) = e^{Fh}$ and $Q(h) = \int_{0}^h e^{F s}LL^\top e^{F^\top s}\dd s$.
    We are interested in computing the posterior distributions $Y(t_k) \mid [\mathcal{Z}(s_j) = 0, j=0, \ldots, k]$, for $k=0, \ldots, K$.
    The resulting state-space model is described by the following transition and likelihood on the measurement process~\eqref{eq:misalignement_measure},
    \begin{equation*}
        \begin{split}
            Y(t_{k+1})&\mid Y(t_k) \sim \calN(A(h)Y(t_k), Q(h)),\\
            \mathcal{Z}(t_k)&\mid Y(t_k) \sim \calN(H_1 Y(t_k)-f(H_0 Y(t_k), t_k), R(t_k)),\\
        \end{split}
    \end{equation*}
    for $k=0,\ldots, K$, $H_{0} = (1, 0, 0, \ldots, 0)$, $H_{1} = (0, 1, 0, \ldots, 0)$, and where a measurement variance $R$ has been added to the measurement process for greater generality.
    However, conditioning $Y(t_k)$ on $\mathcal{Z}(t_j) = 0$ for $j\leq k$ is in general not feasible because of the non-linearity dependency of the measurement process $\mathcal{Z}$ on $Y$ and we need to resort to approximations.
    Different such approximation schemes exist, either Gaussian-based or particle filter-based~\citep[see, e.g.,][]{ChopinSMC2020,sarkka2023bayesian}.
    Here, we focus on the extended Kalman filter which provides closed-form Gaussian approximations to the solution.

    At time $t_k$, let $Y(t_k) \mid [\mathcal{Z}(t_j) = 0, j=0, \ldots, k] \sim \mathcal{N}(m(t_k), P(t_k))$, be given by its approximated filtering posterior.
    At time $t_{k+1}$, we have $Y(t_{k+1}) \mid [\mathcal{Z}(t_j) = 0, j=0, \ldots, k] \sim \mathcal{N}(m^{-}(t_{k+1}), P^{-}(t_{k+1}))$ with
    \begin{equation}
        \label{eq:pred_eq_kalman}
        m^{-}(t_{k+1})\coloneqq A(h)m(t_k),\quad P^{-}(t_{k+1})\coloneqq A(h)P(t_k)A(h)^\top + Q(h).
    \end{equation}
    The updated posterior $Y(t_{k+1})\mid [\mathcal{Z}(t_j) = 0, j=0,\ldots,k+1 ]\sim \calN(m(t_{k+1}), P(t_{k+1}))$ can then be approximated by
    \begin{equation}
        \label{eq:kalman_filter_update_eq}
        \begin{split}
            S(t_{k+1}) &\coloneqq \tilde{H}(t_{k+1})P^{-}(t_{k+1})\tilde{H}(t_{k+1})^\top + R(t_{k+1}),\\
            K(t_{k+1}) &\coloneqq P^{-}(t_{k+1})\tilde{H}(t_{k+1})^{\top} S(t_{k+1})^{-1},\\
            \hat{z}(t_{k+1}) &\coloneq \tilde{H}m^{-}(t_{k+1})-f(H_0 m^{-}(t_{k+1}), t_{k+1}),\\
            m(t_{k+1}) &= m^{-}(t_{k+1}) - K(t_{k+1})\hat{z}(t_{k+1}),\\
            P(t_{k+1}) &= (I_{q+1} - K(t_{k+1})\tilde{H}(t_{k+1}))P^{-}(t_{k+1}).
        \end{split}
    \end{equation}
    Taking $\tilde{H}=H_{1}$ is the extended Kalman filter of order 0 (EKF0) scheme, while the EKF1 is given by $\tilde{H}:t\mapsto H_{1}-J_{f(., t)}(H_{0}m^{-}(t))H_{0}$ where $J_{f(., t)}$ denotes the Jacobian of $f(., t)$~\citep[for a complete review of Gaussian state-space model inference, see, e.g.,][Ch. 5]{sarkka2023bayesian}.

    \subsection{Approximation of the Brownian Motion for Solving SDEs}
    \label{section:brownian_approximation}
    As mentioned in Section~\ref{section:introduction}, we are interested in numerically solving the following additive-noise SDE:
    \begin{equation}
        \begin{split}
            \label{eq:additive_sde2}
            \dd X_{t} &= \mu(X_{t}, t)\dd t + \sigma(t)\dd B_{t}, \indent t \in [0, T],
        \end{split}
    \end{equation}
    whilst preserving uncertainty information regarding its path distribution.
    To apply the existing probabilistic numerics for ODEs to this problem, we will approximate a sample path of~\eqref{eq:additive_sde2} by a path that solves a succession of random ODEs. This is done through the use of differentiable approximations of the Brownian motion. %

    We adopt the approach of~\citet{Foster2020OptimalPolyBrownian} which uses fixed-degree approximations on an ever-thinner grid and then stitches them together to form a global approximation of the Brownian motion. Such an approximation allows the construction of schemes with a high order of convergence for a low computational budget.
    Formally, if the random function $\beta_{k}$ is a given approximation of the Brownian motion over the interval $[t_{k}, t_{k+1}]$, $\beta(t) = \sum_{k} \mathds{1}_{[t_{k}, t_{k+1})}(t) \beta_{k}(t)$ will be a global piecewise approximation of the Brownian motion over the interval $[0, T]$ constructed using those $\beta_{k}$'s.
    Assuming such an approximation for the Brownian motion, we are then interested in solving the following sequence of random ODEs. For $k\in [0, K-1]$,
    \begin{equation}
        \begin{split}
            \label{eq:random_ode}
            \dd{X_{k}}(t) &= \mu(X_{k}(t), t)\dd t+\sigma(t)\dd{\beta_{k}}(t),\indent \textup{for all } t\in [t_{k}, t_{k+1}),\\
            X_{0}(t_{0}) &= Z,\indent
            X_{k}(t_{k}) = X_{k-1}(t_{k}) \textup{ if } k\geq 1.
        \end{split}
    \end{equation}
    Assuming that approximations $\beta_k(t)$ are differentiable, the above can be written in a more standard form
    \begin{equation}
        \begin{split}
            \label{eq:standard_random_ode}
            \frac{\dd{X_{k}}(t)}{\dd t} &= f_{k}(X_k(t), t),\indent \textup{for all } t\in [t_{k}, t_{k+1}),\\
            X_{0}(t_{0}) &= Z,\indent
            X_{k}(t_{k}) = X_{k-1}(t_{k}) \textup{ if } k\geq 1,
        \end{split}
    \end{equation}
    where we have introduced the $(k+1)$-th random vector field $f_k$
    \begin{equation}
        \label{eq:random_vector_field}
        f_{k}(x, t) = \mu(x, t) + \sigma(t)\frac{\dd \beta_{k}}{\dd t} (t),\indent \forall x\in \bbR^d, t\in[t_{k},t_{k+1}].
    \end{equation}
    Solving the previous sequence of ODEs~\eqref{eq:random_ode} is in general not feasible in closed-form, and, instead, we need to resort to a numerical ODE solver to compute $\wtX_{t_1}\approx X_{0}(t_1), \ldots, \allowbreak\wtX_{t_K}\approx X_{K-1}(t_K)$.
    This is summarised in Algorithm~\ref{alg:piecewise-det-solver}.

\begin{figure}[t]
        \centering
        \begin{tikzpicture}[scale=0.7, font=\LARGE]
    \begin{axis}
        [
        legend cell align={left},
        width=.8\textwidth,
        height=.35\textheight,
        xlabel={$t$},
        legend style={
            nodes={scale=0.8, transform shape},
            at={(0,1)},
            anchor=north west,
            fill opacity=0.8, 
            text opacity =1},
        grid=both,
        xmin=0, xmax=0.333,
        ymin=-1.0, ymax=0.7,
        xticklabel style={
        /pgf/number format/fixed,
        /pgf/number format/precision=2
        },
        ]

        \pgfplotstableread[x=t, y=path, col sep=comma]{numerics/path/brownian_fine.csv}\loadedtableI
        \addplot[black, line width=1pt] table {\loadedtableI};
        \pgfplotstableread[x=t, y=path, col sep=comma]{numerics/path/brownian.csv}\loadedtableII
        \addplot[black, dashed, line width=1pt] table {\loadedtableII};
        \pgfplotstableread[x=t, y=path, col sep=comma]{numerics/path/para_traj.csv}\loadedtableIII
        \addplot[gray, line width=1pt] table {\loadedtableIII};
        
        \legend{Fine path, Linear approximation, Parabola approximation};
        
    \end{axis}
\end{tikzpicture}
        \caption{A Brownian path and its linear and parabola approximations. The piecewise derivative of the parabola approximation is discontinuous at the $t_k$'s.}
        \label{fig:Fig0}
    \end{figure}
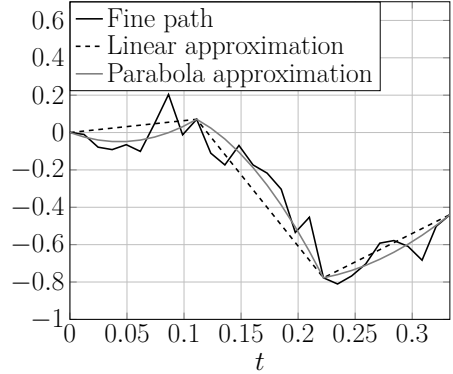

    \begin{algorithm}[t]
        \DontPrintSemicolon
        \caption{A piecewise approximation for simulating a sample path}\label{alg:piecewise-det-solver}
        \SetKwInOut{Input}{input}
        \SetKwInOut{Output}{output}
        \

        \Input{An integration grid $0 = t_{0}<t_{1}<\cdots<t_{K} = T$}
        \Output{An approximate discrete sample path $\wtX_{t_{0}}, \wtX_{t_{1}}, \ldots, \wtX_{t_{K}}$}
        Set $\wtX_{t_{0}} = Z$\;
        \For{$k=0$ \KwTo $K-1$}{
            Sample an approximation $\beta_{k}$ for $t \in [t_k, t_{k+1})$\;
            Compute $\wtX_{t_{k+1}}$ using an ODE solver for the $(k+1)$-th IVP~\eqref{eq:standard_random_ode}\;
        }
        \Return{$\wtX_{t_{0}},\ldots,\wtX_{t_K}$}
    \end{algorithm}
    In practice, the choice of the $\beta_{k}$'s greatly impacts the numerical efficiency of Algorithm~\ref{alg:piecewise-det-solver}, and different choices may lead to different convergence properties or lack thereof.
    \citet{Foster2020OptimalPolyBrownian} introduces a family of $L^2$-optimal polynomial approximation of the Brownian motion that exhibits strong convergence properties: under smoothness assumption on the drift and diffusion functions,
    the solution of~\eqref{eq:random_ode} using the second-degree polynomial approximation converges in $\calL^2$, globally, both strongly and weakly with an order of at least $1.0$ and strongly locally with an order of at least $2.0$~\citep[][Th. 3.14., Th. 3.17]{Foster2020OptimalPolyBrownian}.
    The parabola approximation is, component-wise, given by
    \begin{align}
        \label{eq:parabola}
        \beta_{k}(t) = B_{t_{k}} + B_{t_{k}, t_{k+1}}u+\sqrt{6}u(u-1)I_{t_{k},t_{k+1}},
    \end{align}
    where $u=\frac{t-t_{k}}{t_{k+1}-t_{k}}$ and $B_{t_{k}}, B_{t_{k}, t_{k+1}}, I_{t_{k}, t_{k+1}}$ are independent zero-mean Gaussian random variables with variances $t_{k}$, $t_{k+1}-t_{k}$, and $\frac{t_{k+1}-t_{k}}{2}$, respectively. An illustration of the resulting approximation is given in Figure~\ref{fig:Fig0}.

    \section{Modelling Uncertainty in SDEs}
    \label{section:method}
    In this section, we introduce three novel probabilistic numerics solutions to the additive-noise SDE:
    \begin{equation}
        \begin{split}
            \label{eq:additive_sde3}
            \dd X_{t} &= \mu(X_{t}, t)\dd t + \sigma(t)\dd B_{t}, \indent t \in [0, T]\\
            X_{0} &= Z,
        \end{split}
    \end{equation}
    where $T>0$ is a finite deterministic horizon, $\mu$ and $\sigma$ are vector-valued functions satisfying some latter-specified assumptions, $B$ is a $m$-dimensional Brownian motion, and $Z$ is a square-integrable random variable independent of the filtration induced by $B$.
    The first one, the Gaussian SDE filter, corresponds to applying a single Kalman filter step for each ODE and sampling from the posterior distribution.
    The second one, the Gaussian Mixture SDE filter, is a slightly modified version of the former where, instead of sampling at each time step, we carry over the mean and the variances of the sample path solution. This then provides the marginal distribution of the solution as a Gaussian with random mean and covariance, justifying the name.
    The third one, the Marginalised Gaussian SDE filter, incorporates the random coefficients defining the Brownian approximation as part of the latent state modelling the ODE solution. The resulting inference process then forms a joint distribution over the solution and the driving Brownian path, which we can then discard to obtain the marginal ODE solution.
    As discussed in Section~\ref{section:brownian_approximation}, we will consider a, not yet specified, piecewise random ODE given by~\eqref{eq:random_ode} or equivalently \eqref{eq:standard_random_ode}.

    Importantly, we give strong convergence guarantees for the two first methods when EKF0~\citep{HennigKerstingTronarpSarkka2019} is used on top of the parabola-ODE method~\citep{Foster2020OptimalPolyBrownian}.

    \subsection{The Gaussian SDE Filters}
    For the sake of simplicity, we solely focus on the unidimensional problem, $d=1$, and consider a regular integration grid, $t_{k+1}-t_{k} = \delta$.

    Given a realisation of a Brownian approximation over $[t_k,t_{k+1}]$, our method consists in applying a Gaussian ODE filter to the resulting ODE with vector field a realisation of $f_k$ given by~\eqref{eq:random_vector_field}. We can then sample from the (Gaussian) approximate solution of the ODE at time $t_{k+1}$ and repeat the procedure between $[t_{k+1}, t_{k+2}]$.
    Formally, we introduce, for $t\in [t_k, t_{k+1}]$, the variable $Y_k(t)=\allowbreak(Y_k^{(0)}(t), Y_k^{(1)}(t))^{\top}$ with the initialisation $Y_k(t_k)$ to be sampled from the distribution of $Y^{(0)}_{k-1}(t_k)$:
    \begin{equation*}
        \wtX_{t_k}\sim \mathcal{L}(Y_{k-1}^{(0)}(t_k)),\quad
    \end{equation*}
    where $\mathcal{L}(Y_{k-1}^{(0)}(t_k))$ is the (Gaussian) law of $Y_{k-1}^{(0)}(t_k)$.
    We define the following state-space model,
    \begin{equation*}
        \begin{split}
            Y_k(t_k) &\coloneqq \begin{bmatrix}
                                    \wtX_{k} \\
                                    f_k(\wtX_{k}, t_k)
            \end{bmatrix},\\
            Y_k(t_{k+1})\mid Y_k(t_k) &\sim \calN(A(\delta)Y_{k}(t_k), Q(\delta)),\\
            \mathcal{Z}_k(t_{k+1})\mid Y_k(t_{k+1}) &\sim \calN(Y_k^{(1)}(t_{k+1})-f_k(Y_k^{(0)}(t_{k+1}),t_{k+1}), R),\\
        \end{split}
    \end{equation*}
    where the observation noise $R$ is assumed to be $0$ (see Remarks 1.1, 1.2 and 1.3 in the supplementary material for a generalisation to $R=O(\delta^2)$). The next value of $\wtX_{t_k}$ is then sampled from the distribution of $Y_k^{(0)}(t_{k+1}) \mid \{\mathcal{Z}_k(t_{k+1}) = 0\}$ using one of the Gaussian approximations described in Section~\ref{section:gaussian_ode_filters}, the procedure is then repeated over, resulting in Algorithm~\ref{alg:ekf0-scheme2}.
    \begin{remark}
        The previous approach is different from solving the piecewise ODE on the whole interval $[0, T]$ using~\ref{section:gaussian_ode_filters}.
        Indeed, the discontinuity (see Figure~\ref{fig:Fig0}) at discretisation points of the piecewise derivative of the parabola approximation $\sum_{k=0}^{K-1}\mathds{1}_{[t_k, t_{k+1})}(t)\frac{\dd \beta_k}{\dd s }(t)$ makes the global vector field given by
        \begin{equation*}
            f(x,t)=\mu(x,t)+\sigma(t)\sum_{k=0}^{K-1}\mathds{1}_{[t_k, t_{k+1})}(t)\frac{\dd \beta_k}{\dd s }(t),
        \end{equation*}
        discontinuous at the $t_k$'s.
        This sudden change of the vector field at time $t_k$ constrains us to include in the scheme, both the left and right values, $f_{k-1}(\cdot, t_k)$ and $f_k(\cdot, t_k)$, respectively.
    \end{remark}
    \begin{algorithm}[t]
        \DontPrintSemicolon
        \caption{Gaussian SDE Filter}\label{alg:ekf0-scheme2}
        \SetKwInOut{Input}{input}
        \SetKwInOut{Output}{output}
        \

        \Input{A regular integration grid given by $\delta$\;\\
        }
        \Output{An approximate discrete sample path $\wtX_{t_{0}}, \wtX_{t_{1}}, \ldots, \wtX_{t_{K}}$}
        Set $\wtX_{t_{0}} = Z$\;
        \For{$k=0$ \KwTo $K-1$}{
            Sample an approximation $\beta_{k}$ for $t\in[t_k, t_{k+1})$\;
            Set initial distribution $Y_{k}(t_k)\coloneqq \begin{bmatrix}
                                                              \wtX_{t_{k}}\\ f_{k}(\wtX_{t_{k}},t_{k})
            \end{bmatrix}$\;
            Compute the posterior distribution $Y_{k}(t_{k+1})\sim \calN(m_k(t_{k+1}), P_k(t_{k+1}))$, using one-step EKF0/1\;
            Sample $\wtX_{t_{k+1}}$ from the first component of $\calN(m_k(t_{k+1}), P_k(t_{k+1}))$\;
        }
        \Return $\wtX_{t_0},\ldots, \wtX_{t_{K-1}}$
    \end{algorithm}

    \begin{assumption}[Regularity of the drift and diffusion functions]
        \label{assumption:reg_of_sde_flow}
        The drift function
        \begin{enumerate}
        [label=(\roman*)]
            \item is uniformly bounded with respect to both $x$ and $t$, i.e., for all $x\in \bbR$, $t\in [0, T]$, $\abs{\mu(x,t)}\leq L$,
            \item is continuously differentiable with respect to both $x$ and $t$,
            \item admits uniformly bounded partial derivatives with respect to $x$ and $t$, i.e., for all $x\in \bbR$ and $t\in [0, T]$, $\max\{\abs{\partial_x \mu(x,t)},\abs{\partial_t \mu(x, t)}\} \leq L$,
        \end{enumerate}
        where $L$ is a suitable finite constant.
        Moreover, the diffusion function is continuously differentiable with respect to time.
    \end{assumption}
    Provided the chosen piecewise polynomial approximation is the parabola approximation given by~\eqref{eq:parabola}, and under the Assumption~\ref{assumption:reg_of_sde_flow}, we can prove that Algorithm~\ref{alg:ekf0-scheme2} converges as $\delta \to 0$.
    \begin{theorem}
    [Strong convergence of the Algorithm~\ref{alg:ekf0-scheme2} using EKF0 linearisation]
        \label{th:ekf0-scheme2}
        Assume the drift and diffusion functions satisfy Assumption~\ref{assumption:reg_of_sde_flow}, moreover, assume the Brownian approximation $\beta$ is given by the piecewise parabola~\eqref{eq:parabola}.
        Let $X$ denote the solution of the SDE~\eqref{eq:additive_sde3}.
        Let $\wtX$ denote the solution given by Algorithm~\ref{alg:ekf0-scheme2} using EKF0 linearisation
        with an integrated Ornstein-Uhlenbeck process of parameter $\theta \geq 0$ (see Section~\ref{section:gaussian_ode_filters}) as a prior.
        Then for all $n \in [0, K-1]$:
        \begin{equation*}
            \begin{split}
                \norm{X_{t_1}-\wtX_{t_1}}_{\calL^2} = O(\delta^{1.5}),\quad
                \norm{X_{t_{n+1}}-\wtX_{t_{n+1}}}_{\calL^2} = O(\delta),
            \end{split}
        \end{equation*}
        where $\norm{\cdot}_{\calL^2}=\sqrt{\bbE[\norm{\cdot}^2]}$ is the $\calL^2$-norm and the expectation is taken over the SDE solution's distribution.
    \end{theorem}
    The proof is given in the supplementary material \textbf{Proof 1}.

    \begin{remark}
        While the bounded derivative assumption in Assumption~\ref{assumption:reg_of_sde_flow} is standard, the bounded-ness of the drift is a strong assumption, which we use to control terms of the form $\mu \partial_x \mu$. We however believe that this could be weakened. Indeed, numerical experiments in~\ref{section:experiments} are conducted on a model for which it is not verified, but where the scheme still converges empirically.
    \end{remark}

    \subsection{A Gaussian Mixture Solution}
    \label{section:extension}
    A slightly modified version of the previous scheme, perhaps more in spirit with Gaussian ODE filters, consists in keeping track of the posterior distributions $Y_0(t_1)\sim \calN(m_0(t_1), P_0(t_1))$, $\ldots$ , $Y_{K-1}\sim \calN(m_{K-1}(t_k), P_{K-1}(t_k))$ instead of sampling from it at each iteration.
    To this end, we enforce the posterior mean and covariance of the solution to the $k$-th IVP at time $t_k$, that is $m_{k-1}^{(0)}(t_k)$ and $P_{k-1}(t_k)_{0,0}$, to be the initial mean and covariance of the prior solution to the $(k+1)$-th IVP at time $t_k$, $m_k^{(0)}(t_k)$ and $P_{k}(t_k)_{0,0}$, i.e.,
    \begin{equation*}
        \begin{aligned}[t]
            m_{0}^{(0)}(0) &= Z,\\
            P_{0}(t_0)_{0,0}&\coloneqq 0,
        \end{aligned}
        \qquad
        \begin{aligned}[t]
            m_{k}^{(0)}(t_k) &\coloneqq m_{k-1}^{(0)}(t_k),\\
            P_{k}(t_k)_{0,0} &\coloneqq P_{k-1}(t_k)_{0,0}.
        \end{aligned}
    \end{equation*}
    The vector field components for the mean and covariance are then updated accordingly to the sampled Brownian path on $[t_k, t_{k+1}]$. For EKF1, the approximations are obtained using first-order Taylor series expansion~\citep[see, e.g.,][Ch. 5]{sarkka2023bayesian}:
    \begin{equation}
        \label{eq:update_mean_var_ekf1_scheme2}
        \begin{split}
            m_{k}^{(1)}(t_k)&=\bbE[f_k(Y^{(0)}_k(t_k), t_k)]\\
            &\approx f_k(m_{k}^{(0)}, t_k),\\
            P_{k}(t_k)_{0,1}&=\bbC[Y_k^{(0)}(t_k), f_k(Y^{(0)}(t_k), t_k)]\\
            &\approx P_{k}(t_k)_{0,0}\partial_x f_k(m_k^{(0)}(t_k),t_k),\\
            P_{k}(t_k)_{1,1}&=\bbV[f_k(Y^{(0)}(t_k), t_k)]\\
            &\approx P_{k}(t_k)_{0,0}(\partial_x f_k(m_k^{(0)}(t_k), t_k))^2,
        \end{split}
    \end{equation}
    and for EKF0,
    \begin{equation}
        \label{eq:update_mean_var_ekf0_scheme2}
        \begin{split}
            m_{k}^{(1)}(t_k)\approx f_k(m_{k}^{(0)}, t_k), \quad P_{k}(t_k)_{0,1} \approx 0, \quad P_{k}(t_k)_{1, 1} \approx 0,
        \end{split}
    \end{equation}
    which essentially does not propagate the uncertainty on the solution to its vector field.
    This corresponds to Algorithm~\ref{alg:ekf0-scheme3}.
    \begin{algorithm}[t]
        \DontPrintSemicolon
        \caption{Gaussian Mixture SDE Filter}\label{alg:ekf0-scheme3}
        \SetKwInOut{Input}{input}
        \SetKwInOut{Output}{output}
        \Input{A regular integration grid given by $\delta$\;}
        \Output{A mixture of Gaussian for approximating a discrete sample path $Y_0(t_1), \ldots, Y_{K-1}(t_{K})$}
        \For{$k=0$ \KwTo $K-1$}{
            Sample an approximation $\beta_{k}$ for $t\in [t_k,t_{k+1})$\;
            Set the initial distribution $Y_{k}(t_k)$ to $\calN\left(m_k(t_k), P_k(t_k)\right)$ using~\eqref{eq:update_mean_var_ekf1_scheme2} or~\eqref{eq:update_mean_var_ekf0_scheme2}\;
            Compute the posterior distribution $Y_{k}(t_{k+1})\sim \calN(m_k(t_{k+1}), P_k(t_{k+1}))$, using EKF0/1\;
        }
        \Return $Y_{0}(t_1), \ldots, Y_{K-1}(t_K)$
    \end{algorithm}
    \begin{theorem}
        \label{th:ekf0-scheme3}
        [Strong convergence of Algorithm~\ref{alg:ekf0-scheme3} using EKF0 linearisation]
        Under the Assumptions of Theorem~\ref{th:ekf0-scheme2}, the same convergence orders as~\ref{th:ekf0-scheme2} hold for Algorithm~\ref{alg:ekf0-scheme3} using EKF0 linearisation.
    \end{theorem}
    The proof is given in the supplementary material \textbf{Proof 2}.

    \subsection{Marginalising the Brownian Approximation}
    \label{section:extension2}
    In the previous sections, we successively sampled, for $k=0,\ldots, K-1$, an approximation of the Brownian motion, and then, conditionally on the approximation, the corresponding probabilistic solution $X_k(t_{k+1})$ to the random ODE. This two-step approach, although consistent, is wasteful in that it (i) requires sampling several random variables that are then discarded, (ii) does not provide a closed-form solution for the (approximate) transition density. %
    In this section, we take a different approach, which consists in marginalising the state distribution with respect to the ``random parabola'' explicitly rather than under samples, thus avoiding the sampling procedure of the Brownian approximation required by Algorithms~\ref{alg:ekf0-scheme2} and~\ref{alg:ekf0-scheme3}.
    To do so, we now explicitly assume that the Brownian approximation is given by the parabola approximation of~\citep{Foster2020OptimalPolyBrownian} which we reviewed in Section~\ref{section:brownian_approximation}.

    Marginalising the Brownian approximation then more simply corresponds to computing the joint posterior distribution over the solution and the Gaussian parabola's coefficients $(B_{k\delta, ({k+1})\delta}, I_{k\delta, ({k+1})\delta})$ appearing in~\eqref{eq:parabola}.
    To this end, let us artificially include the Brownian approximation $(B_{k\delta, ({k+1})\delta}, I_{k\delta, ({k+1})\delta})$ inside the state variable $Y_k$,
    \begin{equation*}
        Y_k(t) =
        \begin{pmatrix}
            Y_k^{(0)}(t)               \\
            Y_k^{(1)}(t)               \\
            B_{k\delta, ({k+1})\delta} \\
            I_{k\delta, ({k+1})\delta}
        \end{pmatrix},\quad \forall t\in [t_k,t_{k+1}].
    \end{equation*}
    Similarly to the first scheme, we initialise the first two components of $Y_k$ using a previously sampled solution $\wtX_{t_k}$, for $k\geq 1$:
    \begin{equation*}
        \begin{split}
            Y_k^{(0)}(t_k) \coloneqq \wtX_{t_k}\sim Y_{k-1}^{(0)}(t_k),\quad
            Y_k^{(1)}(t_k) \coloneqq f_k(\wtX_{t_k}, t_k).
        \end{split}
    \end{equation*}
    We construct a specific Kalman filtering step by exploiting the structure of the state-space model with respect to $(B_{k\delta, ({k+1})\delta}, I_{k\delta, ({k+1})\delta})$ to obtain (see Section 1.4 in the supplementary material)
    \begin{equation}
        \label{eq:ssm_3}
        \begin{split}
            Y_k(t_k) &\sim \calN\left(\begin{bmatrix}
                                          \wtX_{t_k}           \\
                                          \mu(\wtX_{t_k}, t_k) \\
                                          0                    \\
                                          0
            \end{bmatrix},
            \begin{bmatrix}
                0 & 0                              & 0           & 0                              \\
                0 & \frac{4}{\delta} \sigma(t_k)^2 & \sigma(t_k) & -\frac{\sqrt{6}}{2}\sigma(t_k) \\
                0 & \sigma(t_k)                    & \delta      & 0                              \\
                0 & -\frac{\sqrt{6}}{2}\sigma(t_k) & 0           & \frac{\delta}{2}
            \end{bmatrix}
            \right),\\
            Y_k(t_{k+1})\mid Y_k(t_k) &\sim \calN(\bar{A}(\delta)Y_{k}(t_k), \bar{Q}(\delta)),\\
            \mathcal{Z}_k(t_{k+1})\mid Y_k(t_{k+1}) &\sim \calN(Y_k^{(1)}(t_{k+1})-\bar{f}_k(Y_k^{(0)}(t_{k+1}),t_{k+1}), R=0),
        \end{split}
    \end{equation}
    where $\bar{A}$ and $\bar{Q}$ are given by
    \begin{equation*}
        \bar{A} = \textup{Diag}(A, 1, 1),\, \bar{Q} = \textup{Diag}(Q, 0, 0),
    \end{equation*}
    $A$ and $Q$ are defined by~\eqref{eq:transition_matrix_ekf}, and $\bar{f}_k$ is the $(k+1)$-th vector field defined using parabola's coefficients as part of the state, i.e.,
    \begin{equation*}
        \bar{f}_k(Y, t) = \mu(Y^{(0)}, t) + \sigma(t)\left(Y^{(2)} + \frac{\sqrt{6}}{\delta} \left(\frac{2(t-t_k)}{\delta} - 1\right)Y^{(3)}\right),
    \end{equation*}
    where $Y^{(2)}$ and $Y^{(3)}$ are the Gaussian coefficients for the parabola approximation.
    The EKF0 scheme is then obtained using~\eqref{eq:kalman_filter_update_eq} with
    \begin{equation*}
        \tilde{H}_k(t) = H_{1}-\left(0, 0, \frac{\sigma(t)}{\delta}, \sigma(t)\frac{\sqrt{6}}{\delta} \left(\frac{2(t-t_k)}{\delta}-1\right)\right),
    \end{equation*}
    while for the EKF1 scheme, we have
    \begin{equation*}
        \tilde{H}_k(t) = H_1- \left(\partial_x\mu(H_{0}m^{-}_k(t), t), 0,  \frac{\sigma(t)}{\delta}, \sigma(t)\frac{\sqrt{6}}{\delta}\left(\frac{2(t-t_k)}{\delta} -1\right)\right).
    \end{equation*}
    The resulting procedure is given by Algorithm~\ref{alg:ekf0-marginal-scheme}.
    \begin{algorithm}[t]
        \DontPrintSemicolon
        \caption{Marginalised Gaussian SDE Filter}\label{alg:ekf0-marginal-scheme}
        \SetKwInOut{Input}{input}
        \SetKwInOut{Output}{output}
        \Input{A regular integration grid given by $\delta$}
        \Output{An approximate discrete sample path $\wtX_{t_{0}}, \wtX_{t_{1}}, \ldots, \wtX_{t_{K}}$}
        Set $\wtX_{t_{0}} = Z$\;
        \For{$k=0$ \KwTo $K-1$}{
            Set initial distribution $Y_k(t_k)$ to~\eqref{eq:ssm_3} using $\wtX_{t_k}$\;
            Sample $\wtX_{t_{k+1}}$ from the posterior distribution $Y_k(t_{k+1})\sim \calN(m_k(t_{k+1}), P_k(t_{k+1}))$, computed using EKF0/1\;
        }
        \Return $\wtX_{t_0}, \ldots, \wtX_{t_K}$
    \end{algorithm}
    Contrary to Algorithms~\ref{alg:ekf0-scheme2} and~\ref{alg:ekf0-scheme3}, we did not obtain convergence guarantees for Algorithm~\ref{alg:ekf0-marginal-scheme}.
    This is because incorporating the Brownian path as part of the solution makes its treatment largely more complicated, with the non-linearity impacting additional parts of the algorithm.
    However, we empirically observe in Section~\ref{section:experiments} that it enjoys good weak convergence properties, and, therefore, should be considered as a viable algorithm for uncertainty modelling.
    By the Rao-Blackwell theorem, if one were able to perform the marginalisation step exactly (which we are not given the non-linearity involved, and so the consequence has to be understood heuristically), such moment estimates exhibit lower variance compared to those obtained via simulation, see, e.g.,~\citet{sarkka2023bayesian}[Ch. 7.5].

    \subsection{Uncertainty Calibration}
    \label{section:calibration}
    The uncertainty calibration can be carried out by choosing the prior, in particular, the drift and diffusion matrices~\eqref{eq:transition_matrix_ekf}.
    While we restrict ourselves to the set of IOUP processes for simplicity, the procedure we discuss below on the calibration of the diffusion coefficient $\eta$ can be extended to other priors.
    The most adopted method for uncertainty calibration in the Bayesian ODE solver literature is the one discussed in~\citet{HennigKerstingTronarpSarkka2019}, in which a post-hoc point estimate on the diffusion coefficient $\eta$ is obtained under the assumption that the vector field is affine.
    In the affine SDE setting, this condition is equivalent to assuming the drift function $\mu$ is affine with respect to $x$.
    \begin{proposition}[Quasi maximum-likelihood for the diffusion coefficient]
        \label{prop:ml_post_hoc}
        Assume the drift is an affine function with respect to $X$, i.e.,~\eqref{eq:additive_sde2} becomes
        \begin{equation*}
            \begin{split}
                \mathrm{d}X_t &= (g(t)+f(t)X_t)\mathrm{d}t + \sigma(t)\mathrm{d}B_t,\\
                X_0 &= Z.
            \end{split}
        \end{equation*}
        At time $t_k$, let $m_k(t_{k+1}), P_k(t_{k+1})$ denote the posterior mean and covariance given by Algorithm~\ref{alg:ekf0-scheme3} with transition parameters $(A, Q)$.
        Similarly, let $(m_k'(t_{k+1}), P_k'(t_{k+1}))$ be the output of Algorithm~\ref{alg:ekf0-scheme3} but now with $\eta$-scaled transition parameters $(A, Q'=\eta^2Q)$.
        Then, the posterior means coincide, i.e., $m_k' = m_k$ and the posterior covariances $P_k$ and $P_k'$ are proportional with $P_k' = \eta^2 P_k$.
        Additionally, denote the predicted mean and covariance of the residual error $Z_k(t_{k+1})$ by $\hat{z}_k(t_{k+1})$ and $S_k(t_{k+1})$.
        Then, conditionally on the Brownian path and the initialisation, the maximum-likelihood estimate for $\hat{\eta}^2$ given we observe the measurement error $\mathcal{Z}_{k}(t_{k+1}) = 0$ for $k\in [0, K-1]$ is given by
        \begin{equation}
            \label{eq:mle_estmiate_for_the_diffusion_coefficient}
            \begin{split}
                \hat{\eta}^2 = \frac{1}{K}\sum_{k=0}^{K-1} \hat{z}_{k}(t_{k+1})^{\top}S_k(t_{k+1})^{-1}\hat{z}_{k}(t_{k+1}).
            \end{split}
        \end{equation}
    \end{proposition}
    Other calibration methods including inline calibration allowing time-varying diffusion coefficients are discussed in~\citet{Bosch2020CalibratedAP}.

    The calibration uncertainty of a trajectory then goes as follows: (1) use Algorithm~\ref{alg:ekf0-scheme3} to sample from a trajectory, (2) apply~\eqref{eq:mle_estmiate_for_the_diffusion_coefficient} to calibrate the posterior uncertainty on the solution.
    \begin{remark}
        As outlined in Proposition~\ref{prop:ml_post_hoc}
        this procedure is an exact MLE for affine drifts only, but is biased otherwise, hence the name quasi maximum-likelihood. Nonetheless, it provides a fast, and principled, way to calibrate uncertainties for the trajectory.
    \end{remark}

    \section{Experiments}
    \label{section:experiments}
    Throughout this section, we assess the strong order of convergence for the Algorithms~\ref{alg:ekf0-scheme2} and~\ref{alg:ekf0-scheme3} using EKF0 linearisation (see Theorems~\ref{th:ekf0-scheme2} and~\ref{th:ekf0-scheme3}). We empirically measure the weak convergence performances of Algorithms~\ref{alg:ekf0-scheme2},~\ref{alg:ekf0-scheme3} and~\ref{alg:ekf0-marginal-scheme}.
    In addition, we provide an illustration of the calibration method given in Proposition~\ref{prop:ml_post_hoc}.
    The JAX code to reproduce the experiments listed below can be found at the following address: \url{https://github.com/ylefay/bayesianSDEsolver}.

    \label{section:exp_error_def}

    \subsection{Experimental methodology}
    The fine solution we will use as a comparison to our approximate solution $\wtX^{\delta}$ is obtained using the Euler--Maruyama scheme for a time-step $\delta^2$ and is denoted by $X^{em, \delta^2}$.
    By simulating parabola approximations~\eqref{eq:parabola} with scale $\delta^2$, $(B_{k\delta^2, (k+1)\delta^2}, I_{k\delta^2, (k+1)\delta^2})$, we are able to compute the conditional coefficients $(B_{k\delta}, I_{k\delta})\mid (B_{k\delta^2, (k+1)\delta^2}, I_{k\delta^2, (k+1)\delta^2})$, which will be used to compute $\wtX^{\delta}$ by either Algorithm~\ref{alg:ekf0-scheme2} or Algorithm~\ref{alg:ekf0-scheme3}.
    For the state-space model of Algorithm~\ref{alg:ekf0-marginal-scheme}, we compute only weak error estimates: this is because no strong convergence can be computed as the path has been marginalised.

    The number of sample paths for each experiment is $N = 10^5$ and the considered step sizes are $\delta\in \{2^{-4}, 2^{-5}, \ldots, 2^{-10}\}$.
    The estimates of the strong local and global errors are given by
    \begin{equation*}
        \begin{split}
            \hat{\varepsilon}_g = \frac{1}{N}\sum_{k=1}^{N}\norm{\wtX^{\delta}_{T,k}-X^{em, \delta^{2}}_{T,k}},\indent
            \hat{\varepsilon}_l = \frac{1}{N}\sum_{k=1}^{N}\norm{\wtX^{\delta}_{\delta,k}-X^{em, \delta^{2}}_{\delta,k}}.
        \end{split}
    \end{equation*}
    The weak error is estimated using the polynomial function $g(X) = XX^{\top}$,
    \begin{equation*}
        \hat{\varepsilon}_{wg} = \frac{1}{N}\bignorm{\sum_{k=1}^{N}g(\wtX^{\delta}_{T,k})-g(X^{em, \delta^{2}}_{T,k})}.
    \end{equation*}
    Throughout the experiments, we consider the test model defined by the following SDE
    \begin{align*}
        \dd X_0(t) &= 1/\varepsilon(s+X_0(t)-X_0(t)^3-X_1(t))\dd t,\\
        \dd X_1(t) &= (\alpha+\gamma X_0(t) - X_1(t))\dd t +\sigma \dd B_t,\\
        X(0) &= (0, 0)^{\top}, \qquad \quad \,\,\,\,\,
    \end{align*}
    where $\varepsilon=0.1$, $s=0$, $\gamma=1.5$, $\alpha=0.8$, $\sigma=0.3$, $t\in [0, 1]$.
    This corresponds to the stochastic Fitzhugh–Nagumo model, which is a stochastic oscillator used to describe spiking neurons~\citep{FITZHUGH1961445, Berglund2006}. %
    Thus, the transition density's covariance showcases a high-order $O(\delta^3)$ term on the first component and high-order schemes are needed to capture the variance structure (i.e., a $1.5$ weak-order scheme).

    \subsection{Convergence Rates}
    We now turn to empirically assessing the convergence orders.
    Let $\hat{\varepsilon}$ be either the weak error or strong error estimate.
    In the perfect case where the error estimate is of the form $\hat{\varepsilon}=b\delta^{a}$, which is the case for a convergence with order $\alpha$, we have $\ln{\hat{\varepsilon}} = a\ln\delta + \ln{b}$.
    We therefore assume the following linear model with Gaussian residuals:
    \begin{equation}
        \label{eq:linear_model}
        \ln\bbE[\hat{\varepsilon}] = a\ln\delta + \ln b + \nu,
    \end{equation}
    where $\nu\sim \calN(0, c^2)$ and perform a least-square regression on our data to estimate $a$ and $\ln b$. The estimators $(\hat{a}, \hat{b})$, as well as the estimated standard deviations $\hat{c}$ and the $R^2$ coefficients are given in Tables~\ref{tab:tab1},~\ref{tab:tab2} and~\ref{tab:tab3}.
    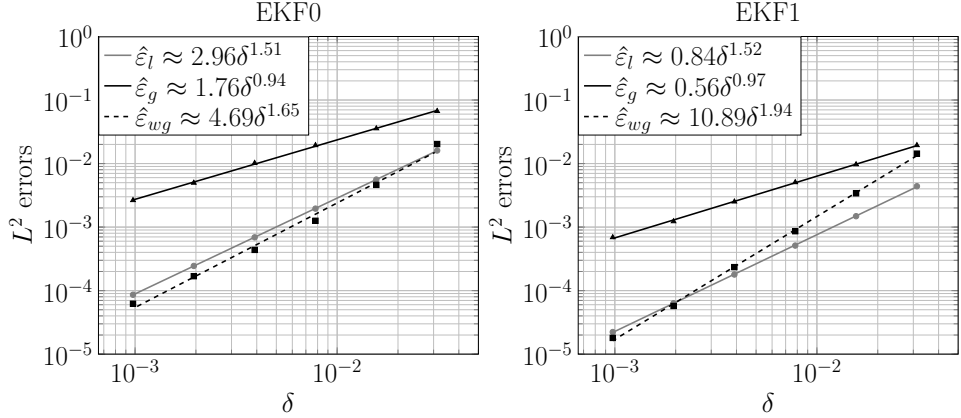
\begin{figure}[t]
        \centering
        \begin{subfigure}[b]{0.47\textwidth}
            \centering
            \begin{tikzpicture}[scale=0.58, font=\LARGE]
    \begin{axis}
        [
        legend cell align={left},
        xlabel={$\delta$},
        ylabel={$L^2$ errors},
        title={EKF0},
        legend style={
            nodes={scale=0.8, transform shape},
            at={(0,1)},
            anchor=north west,
            fill opacity=0.8, 
            text opacity =1},
        grid=both,
        xmode=log,
        ymode=log,
        xmax=0.05,
        ymin=1e-5, ymax=1,
        ]

        \pgfplotstableread[x=deltas, y=errors, col sep=comma]{numerics/EKF0/EKF0_FHN_STRONG_LOCAL_ERRORS.csv}\loadedtableI
        \addplot[gray, only marks, mark=*, mark size=1.7pt] table {\loadedtableI};
        \addplot[gray, line width=1pt]
        table [y={create col/linear regression={y=errors}}] {\loadedtableI};
        \xdef\slopeI{\pgfplotstableregressiona}
        \xdef\bI{\pgfplotstableregressionb}

        \pgfplotstableread[x=deltas, y=errors, col sep=comma]{numerics/EKF0/EKF0_FHN_STRONG_GLOBAL_ERRORS.csv}\loadedtableII
        \addplot[black, only marks, mark=triangle*, mark size=1.7pt] table {\loadedtableII};
        \addplot[black, line width=1pt]
        table [y={create col/linear regression={y=errors}}] {\loadedtableII};
        \xdef\slopeII{\pgfplotstableregressiona}
        \xdef\bII{\pgfplotstableregressionb}

        \pgfplotstableread[x=deltas, y=errors, col sep=comma]{numerics/EKF0/EKF0_FHN_WEAK_GLOBAL_ERRORS.csv}\loadedtableIII
        \addplot[black, only marks, mark=square*, mark size=1.7pt] table {\loadedtableIII};
        \addplot[black, dashed, line width=1pt]
        table [y={create col/linear regression={y=errors}}] {\loadedtableIII};
        \xdef\slopeIII{\pgfplotstableregressiona}
        \xdef\bIII{\pgfplotstableregressionb}

        \legend{,$\hat{\varepsilon}_{l}\approx\pgfmathparse{exp(\bI)} \xdef\ebI{\pgfmathresult} \pgfmathprintnumber[precision=2]{\ebI}\delta^{\pgfmathprintnumber[precision=2, fixed zerofill]
        {\slopeI}}$,,$\hat{\varepsilon}_{g}\approx\pgfmathparse{exp(\bII)} \xdef\ebII{\pgfmathresult} \pgfmathprintnumber[precision=2]{\ebII}\delta^{\pgfmathprintnumber[precision=2, fixed zerofill]
        {\slopeII}}$,,$\hat{\varepsilon}_{wg}\approx\pgfmathparse{exp(\bIII)} \xdef\ebIII{\pgfmathresult} \pgfmathprintnumber[precision=2]{\ebIII}\delta^{\pgfmathprintnumber[precision=2, fixed zerofill]
        {\slopeIII}}$};

    \end{axis}
\end{tikzpicture}
        \end{subfigure}
        \begin{subfigure}[b]{0.47\textwidth}
            \centering
            \begin{tikzpicture}[scale=0.58, font=\LARGE]
    \begin{axis}
        [%
        legend cell align={left},
        xlabel={$\delta$},
        ylabel={$L^2$ errors},
        title={EKF1},
        legend style={
            nodes={scale=0.8, transform shape},
            at={(0,1)},
            anchor=north west,
            fill opacity=0.8, 
            text opacity =1},
        grid=both,
        xmode=log,
        ymode=log,
        xmax=0.05,
        ymin=1e-5, ymax=1,
        ]
        \pgfplotstableread[x=deltas, y=errors, col sep=comma]{numerics/EKF1/EKF1_FHN_STRONG_LOCAL_ERRORS.csv}\loadedtableI
        \addplot[gray, only marks, mark=*, mark size=1.7pt] table {\loadedtableI};
        \addplot[gray, line width=1pt] table [y={create col/linear regression={y=errors}}] {\loadedtableI};
        \xdef\slopeI{\pgfplotstableregressiona}
        \xdef\bI{\pgfplotstableregressionb}
        \pgfplotstableread[x=deltas, y=errors, col sep=comma]{numerics/EKF1/EKF1_FHN_STRONG_GLOBAL_ERRORS.csv}\loadedtableII
        \addplot[black, only marks, mark=triangle*, mark size=1.7pt] table {\loadedtableII};
        \addplot[black, line width=1pt]
        table [y={create col/linear regression={y=errors}}] {\loadedtableII};
        \xdef\slopeII{\pgfplotstableregressiona}
        \xdef\bII{\pgfplotstableregressionb}
        \pgfplotstableread[x=deltas, y=errors, col sep=comma]{numerics/EKF1/EKF1_FHN_WEAK_GLOBAL_ERRORS.csv}\loadedtableIII
        \addplot[black, only marks, mark=square*, mark size=1.7pt] table {\loadedtableIII};
        \addplot[black, dashed, line width=1pt]
        table [y={create col/linear regression={y=errors}}] {\loadedtableIII};
        \xdef\slopeIII{\pgfplotstableregressiona}
        \xdef\bIII{\pgfplotstableregressionb}
        \legend{,$\hat{\varepsilon}_{l}\approx\pgfmathparse{exp(\bI)} \xdef\ebI{\pgfmathresult} \pgfmathprintnumber[precision=2]{\ebI}\delta^{\pgfmathprintnumber[precision=2, fixed zerofill]
        {\slopeI}}$,,$\hat{\varepsilon}_{g}\approx\pgfmathparse{exp(\bII)} \xdef\ebII{\pgfmathresult} \pgfmathprintnumber[precision=2]{\ebII}\delta^{\pgfmathprintnumber[precision=2, fixed zerofill]
        {\slopeII}}$,,$\hat{\varepsilon}_{wg}\approx\pgfmathparse{exp(\bIII)} \xdef\ebIII{\pgfmathresult} \pgfmathprintnumber[precision=2]{\ebIII}\delta^{\pgfmathprintnumber[precision=2, fixed zerofill]
        {\slopeIII}}$};
    \end{axis}
\end{tikzpicture}
        \end{subfigure}
        \caption{Algorithm~\ref{alg:ekf0-scheme2}, Fitzhugh--Nagumo model. Independent IBM prior on each coordinates with $\eta=1$. Strong local and global errors, and weak global errors (grey points, black triangles and black squares, respectively) with log-log regression (grey, black line and dashed black line, respectively).}
        \label{fig:Fig1}
    \end{figure}
    \begin{table}[t]
        \centering
        \begin{tabular}{|c|c|c|c|c|}
            \hline
            $\hat{\varepsilon}$                     & $\hat{a}$ & $\hat{b}$ & $\hat{c}$ & $R^2$   \\
            \hline
            $\hat{\varepsilon}_{l, \textrm{EKF0}}$  & $1.508$   & $1.086$   & $0.007$   & $0.999$ \\
            $\hat{\varepsilon}_{g, \textrm{EKF0}}$  & $0.936$   & $0.565$   & $0.033$   & $0.999$ \\
            $\hat{\varepsilon}_{wg, \textrm{EKF0}}$ & $1.645$   & $1.545$   & $0.216$   & $0.992$ \\
            \hline
            $\hat{\varepsilon}_{l, \textrm{EKF1}}$  & $1.523$   & $-0.171$  & $0.018$   & $0.999$ \\
            $\hat{\varepsilon}_{g, \textrm{EKF1}}$  & $0.972$   & $-0.586$  & $0.036$   & $0.999$ \\
            $\hat{\varepsilon}_{wg, \textrm{EKF1}}$ & $1.935$   & $2.388$   & $0.081$   & $0.999$ \\
            \hline
        \end{tabular}
        \caption{Linear regression summary of the model~\eqref{eq:linear_model} for the Algorithm~\ref{alg:ekf0-scheme2} using either EKF0 or EKF1.}
        \label{tab:tab1}
    \end{table}
    \begin{figure}[t]
        \centering
        \begin{subfigure}[b]{0.47\linewidth}
            \centering
            \begin{tikzpicture}[scale=0.58, font=\LARGE]
    \begin{axis}
        [
        legend cell align={left},
        xlabel={$\delta$},
        ylabel={$L^2$ errors},
        title={EKF0},
        legend style={
            nodes={scale=0.8, transform shape},
            at={(0,1)},
            anchor=north west,
            fill opacity=0.8, 
            text opacity =1},
        grid=both,
        xmode=log,
        ymode=log,
        xmax=0.05,
        ymin=1e-5, ymax=1,
        ]

        \pgfplotstableread[x=deltas, y=errors, col sep=comma]{numerics/EKF0_2/EKF0_2_Attempt2FHN_STRONG_LOCAL_ERRORS.csv}\loadedtableI
        \addplot[gray, only marks, mark=*, mark size=1.7pt] table {\loadedtableI};
        \addplot[gray, line width=1pt]
        table [y={create col/linear regression={y=errors}}] {\loadedtableI};
        \xdef\slopeI{\pgfplotstableregressiona}
        \xdef\bI{\pgfplotstableregressionb}

        \pgfplotstableread[x=deltas, y=errors, col sep=comma]{numerics/EKF0_2/EKF0_2_Attempt2FHN_STRONG_GLOBAL_ERRORS.csv}\loadedtableII
        \addplot[black, only marks, mark=triangle*, mark size=1.7pt] table {\loadedtableII};
        \addplot[black, line width=1pt]
        table [y={create col/linear regression={y=errors}}] {\loadedtableII};
        \xdef\slopeII{\pgfplotstableregressiona}
        \xdef\bII{\pgfplotstableregressionb}

        \pgfplotstableread[x=deltas, y=errors, col sep=comma]{numerics/EKF0_2/EKF0_2_Attempt2FHN_WEAK_GLOBAL_ERRORS.csv}\loadedtableIII
        \addplot[black, only marks, mark=square*, mark size=1.7pt] table {\loadedtableIII};
        \addplot[black, dashed, line width=1pt]
        table [y={create col/linear regression={y=errors}}] {\loadedtableIII};
        \xdef\slopeIII{\pgfplotstableregressiona}
        \xdef\bIII{\pgfplotstableregressionb}

        \legend{,$\hat{\varepsilon}_{l}\approx\pgfmathparse{exp(\bI)} \xdef\ebI{\pgfmathresult} \pgfmathprintnumber[precision=2]{\ebI}\delta^{\pgfmathprintnumber[precision=2, fixed zerofill]
        {\slopeI}}$,,$\hat{\varepsilon}_{g}\approx\pgfmathparse{exp(\bII)} \xdef\ebII{\pgfmathresult} \pgfmathprintnumber[precision=2]{\ebII}\delta^{\pgfmathprintnumber[precision=2, fixed zerofill]
        {\slopeII}}$,,$\hat{\varepsilon}_{wg}\approx\pgfmathparse{exp(\bIII)} \xdef\ebIII{\pgfmathresult} \pgfmathprintnumber[precision=2]{\ebIII}\delta^{\pgfmathprintnumber[precision=2, fixed zerofill]
        {\slopeIII}}$};

    \end{axis}
\end{tikzpicture}
        \end{subfigure}
        \begin{subfigure}[b]{0.47\linewidth}
            \centering
            \begin{tikzpicture}[scale=0.58, font=\LARGE]
    \begin{axis}
        [
        legend cell align={left},
        xlabel={$\delta$},
        ylabel={$L^2$ errors},
        title={EKF1},
        legend style={
            nodes={scale=0.8, transform shape},
            at={(0,1)},
            anchor=north west,
            fill opacity=0.8, 
            text opacity =1},
        grid=both,
        xmode=log,
        ymode=log,
        xmax=0.05,
        ymin=1e-5, ymax=1,
        ]

        \pgfplotstableread[x=deltas, y=errors, col sep=comma]{numerics/EKF1_2/EKF1_2_FHN_STRONG_LOCAL_ERRORS.csv}\loadedtableI
        \addplot[gray, only marks, mark=*, mark size=1.7pt] table {\loadedtableI};
        \addplot[gray, line width=1pt]
        table [y={create col/linear regression={y=errors}}] {\loadedtableI};
        \xdef\slopeI{\pgfplotstableregressiona}
        \xdef\bI{\pgfplotstableregressionb}

        \pgfplotstableread[x=deltas, y=errors, col sep=comma]{numerics/EKF1_2/EKF1_2_FHN_STRONG_GLOBAL_ERRORS.csv}\loadedtableII
        \addplot[black, only marks, mark=triangle*, mark size=1.7pt] table {\loadedtableII};
        \addplot[black, line width=1pt]
        table [y={create col/linear regression={y=errors}}] {\loadedtableII};
        \xdef\slopeII{\pgfplotstableregressiona}
        \xdef\bII{\pgfplotstableregressionb}

        \pgfplotstableread[x=deltas, y=errors, col sep=comma]{numerics/EKF1_2/EKF1_2_FHN_WEAK_GLOBAL_ERRORS.csv}\loadedtableIII
        \addplot[black, only marks, mark=square*, mark size=1.7pt] table {\loadedtableIII};
        \addplot[black, dashed, line width=1pt]
        table [y={create col/linear regression={y=errors}}] {\loadedtableIII};
        \xdef\slopeIII{\pgfplotstableregressiona}
        \xdef\bIII{\pgfplotstableregressionb}

        \legend{,$\hat{\varepsilon}_{l}\approx\pgfmathparse{exp(\bI)} \xdef\ebI{\pgfmathresult} \pgfmathprintnumber[precision=2]{\ebI}\delta^{\pgfmathprintnumber[precision=2, fixed zerofill]
        {\slopeI}}$,,$\hat{\varepsilon}_{g}\approx\pgfmathparse{exp(\bII)} \xdef\ebII{\pgfmathresult} \pgfmathprintnumber[precision=2]{\ebII}\delta^{\pgfmathprintnumber[precision=2, fixed zerofill]
        {\slopeII}}$,,$\hat{\varepsilon}_{wg}\approx\pgfmathparse{exp(\bIII)} \xdef\ebIII{\pgfmathresult} \pgfmathprintnumber[precision=2]{\ebIII}\delta^{\pgfmathprintnumber[precision=2, fixed zerofill]
        {\slopeIII}}$};

    \end{axis}
\end{tikzpicture}
        \end{subfigure}
        \caption{Algorithm~\ref{alg:ekf0-scheme3}, Fitzhugh--Nagumo model. Strong, local and global errors, and weak global errors.}
        \label{fig:Fig2}
    \end{figure}
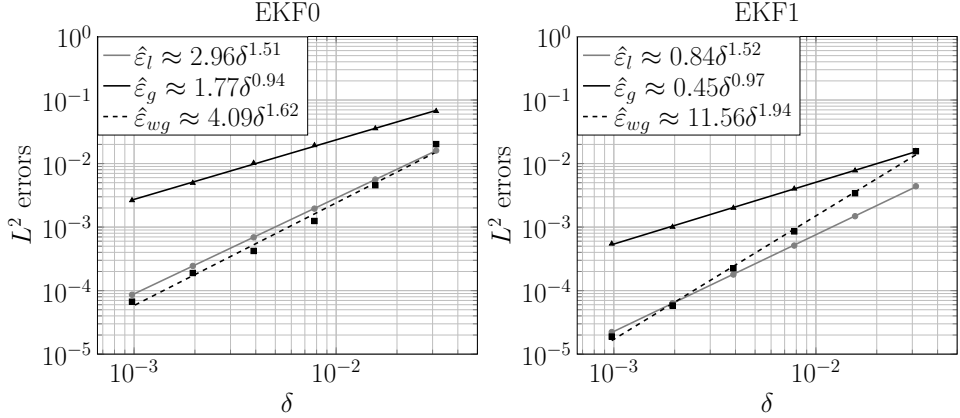
    \begin{table}[t]
        \centering
        \begin{tabular}{|c|c|c|c|c|}
            \hline
            $\hat{\varepsilon}$                     & $\hat{a}$    & $\hat{b}$     & $\hat{c}$    & $R^2$                \\
            \hline
            $\hat{\varepsilon}_{l, \textrm{EKF0}}$  & $1.508$      & $1.086$       & $0.007$      & $0.999$              \\
            $\hat{\varepsilon}_{g, \textrm{EKF0}}$  & $0.937$      & $0.569$       & $0.034$      & $0.999$              \\
            $\hat{\varepsilon}_{wg, \textrm{EKF0}}$ & $1.615$      & $1.409$       & $0.249$      & $0.989$              \\
            \hline
            $\hat{\varepsilon}_{l, \textrm{EKF1}}$  & $1.523$      & $-0.171$      & $0.0176$     & $0.999$              \\
            $\hat{\varepsilon}_{g, \textrm{EKF1}}$  & $0.973$      & $-0.803$      & $0.026$      & $0.999$              \\
            $\hat{\varepsilon}_{wg, \textrm{EKF1}}$ & $1.943$      & $2.447$       & $0.120$      & $0.998$              \\
            \hline
        \end{tabular}
        \caption{Linear regression summary of the model~\eqref{eq:linear_model} for Algorithm~\ref{alg:ekf0-scheme3} using either EKF0 or EKF1.}
        \label{tab:tab2}
    \end{table}
    \begin{figure}[t]
        \centering
        \begin{subfigure}[t]{0.47\textwidth}
            \centering
            \begin{tikzpicture}[scale=0.58, font=\LARGE]
    \begin{axis}
        [
        legend cell align={left},
        xlabel={$\delta$},
        ylabel={$L^2$ errors},
        title={EKF0},
        legend style={
            nodes={scale=0.8, transform shape},
            at={(0,1)},
            anchor=north west,
            fill opacity=0.8, 
            text opacity =1},
        grid=both,
        xmode=log,
        ymode=log,
        xmax=0.05,
        ymin=1e-5, ymax=1,
        ]
        \pgfplotstableread[x=deltas, y=errors, col sep=comma]{numerics/EKF0/EKF0_FHN_WEAK_GLOBAL_ERRORS.csv}\loadedtableI
        \addplot[black, only marks, mark=square*, mark size=1.7pt] table {\loadedtableI};
        \addplot[black, dashed, line width=1pt]
        table [y={create col/linear regression={y=errors}}] {\loadedtableI};
        \xdef\slopeI{\pgfplotstableregressiona}
        \xdef\bI{\pgfplotstableregressionb}

        \pgfplotstableread[x=deltas, y=errors, col sep=comma]{numerics/SSM_EKF0/EKF0_SSM_FHN_WEAK_GLOBAL_ERRORS.csv}\loadedtableIII
        \addplot[gray, only marks, mark=square*, mark size=1.7pt] table {\loadedtableIII};
        \addplot[gray, dashed, line width=1pt]
        table [y={create col/linear regression={y=errors}}] {\loadedtableIII};
        \xdef\slopeIII{\pgfplotstableregressiona}
        \xdef\bIII{\pgfplotstableregressionb}

        \legend{,$\hat{\varepsilon}_{wg,(2)}\approx\pgfmathparse{exp(\bI)} \xdef\ebI{\pgfmathresult} \pgfmathprintnumber[precision=2]{\ebI}\delta^{\pgfmathprintnumber[precision=2, fixed zerofill]
        {\slopeI}}$,,$\hat{\varepsilon}_{wg,(4)}\approx\pgfmathparse{exp(\bIII)} \xdef\ebIII{\pgfmathresult} \pgfmathprintnumber[precision=2]{\ebIII}\delta^{\pgfmathprintnumber[precision=2, fixed zerofill]
        {\slopeIII}}$};

    \end{axis}
\end{tikzpicture}
        \end{subfigure}
        \begin{subfigure}[t]{0.47\textwidth}
            \centering
            \begin{tikzpicture}[scale=0.6, font=\LARGE]
    \begin{axis}
        [
        legend cell align={left},
        xlabel={$t$},
        ylabel={Trajectory estimated means},
        title={EKF1},
        legend style={
            nodes={scale=0.8, transform shape},
            at={(0,1)},
            anchor=north west,
            fill opacity=0.8, 
            text opacity =1},
        grid=both,
        xmin=0, xmax=1,
        ymin=-1, ymax=1,
        ]

        \pgfplotstableread[x=t, y=mean, col sep=comma]{numerics/SSM_EKF1/EKF1_SSM_FHN_correct_mean_first.csv}\loadedtableI
        \addplot[black, line width=1pt] table {\loadedtableI};
        \pgfplotstableread[x=t, y=mean, col sep=comma]{numerics/SSM_EKF1/EKF1_SSM_FHN_correct_mean_secnd.csv}\loadedtableII
        \addplot[darkgray, line width=1pt] table {\loadedtableII};
        \pgfplotstableread[x=t, y=mean, col sep=comma]{numerics/SSM_EKF1/EKF1_SSM_FHN_incorrect_mean_first.csv}\loadedtableIII
        \addplot[gray, line width=1pt] table {\loadedtableIII};
        \pgfplotstableread[x=t, y=mean, col sep=comma]{numerics/SSM_EKF1/EKF1_SSM_FHN_incorrect_mean_secnd.csv}\loadedtableIV
        \addplot[coolgrey, line width=1pt] table {\loadedtableIV};

        \legend{$1^{\textup{st}}$ coord. FS, $2^{\textup{nd}}$ coord. FS, $1^{\textup{st}}$ coord., $2^{\textup{nd}}$ coord.};
        
    \end{axis}
\end{tikzpicture}
        \end{subfigure}
        \caption{Fitzhugh--Nagumo model. Left: weak global errors for Algorithm~\ref{alg:ekf0-scheme2} and~\ref{alg:ekf0-marginal-scheme} using EKF0, in log-log scale. Right: trajectory means, from the fine solution (FS) and the estimated solution given by Algorithm~\ref{alg:ekf0-marginal-scheme} using EKF1 with $\delta=2^{-10}$.}
        \label{fig:Fig3}
    \end{figure}
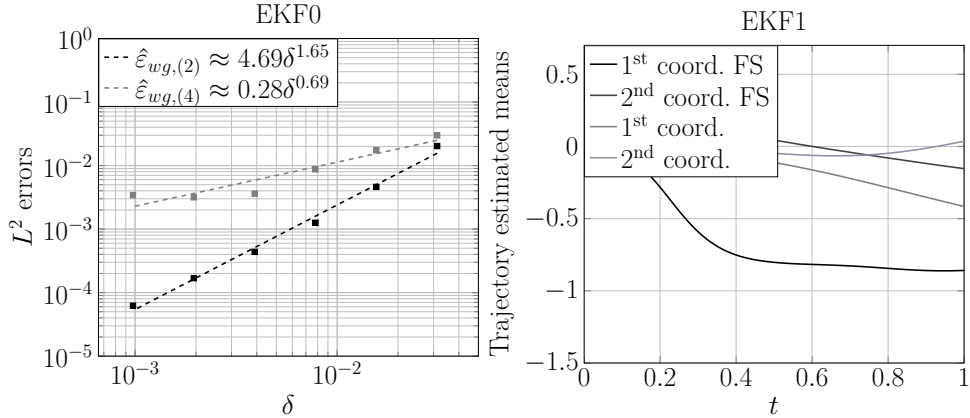
    \begin{table}[t]
        \centering
        \begin{tabular}{|c|c|c|c|c|}
            \hline
            $\hat{\varepsilon}$                     & $\hat{a}$    & $\hat{b}$    & $\hat{c}$    & $R^2$                \\
            \hline
            $\hat{\varepsilon}_{wg, \textrm{EKF0}}$ & $0.694$      & $-1.287$     & $0.350$      & $0.892$              \\
            \hline
        \end{tabular}
        \caption{Linear regression summary of the model~\eqref{eq:linear_model} for Algorithm~\ref{alg:ekf0-marginal-scheme} using EKF0.}
        \label{tab:tab3}
    \end{table}
    \begin{figure}[t]
        \centering
        \begin{subfigure}[b]{0.47\textwidth}
            \centering
            \begin{tikzpicture}[scale=0.58, font=\LARGE]
    \begin{axis}
        [
        legend cell align={left},
        xlabel={$\delta$},
        ylabel={$L^2$ errors},
        title={Strong error comparisons},
        legend style={
            nodes={scale=0.8, transform shape},
            at={(0,1)},
            anchor=north west,
            fill opacity=0.8, 
            text opacity =1},
        grid=both,
        xmode=log,
        ymode=log,
        xmax=0.05,
        ymin=1e-5, ymax=1,
        ]

        \pgfplotstableread[x=deltas, y=errors, col sep=comma]{numerics/EM/EM_FHN_STRONG_GLOBAL_ERRORS.csv}\loadedtableI
        \addplot[black, only marks, mark=triangle*, mark size=1.7pt] table {\loadedtableI};
        \addplot[black, line width=1pt]
        table [y={create col/linear regression={y=errors}}] {\loadedtableI};
        \xdef\slopeI{\pgfplotstableregressiona}
        \xdef\bI{\pgfplotstableregressionb}

        \pgfplotstableread[x=deltas, y=errors, col sep=comma]{numerics/EKF0/EKF0_FHN_STRONG_GLOBAL_ERRORS.csv}\loadedtableII
        \addplot[darkgray, only marks, mark=square*, mark size=1.7pt] table {\loadedtableII};
        \addplot[darkgray, line width=1pt]
        table [y={create col/linear regression={y=errors}}] {\loadedtableII};
        \xdef\slopeII{\pgfplotstableregressiona}
        \xdef\bII{\pgfplotstableregressionb}

        \pgfplotstableread[x=deltas, y=errors, col sep=comma]{numerics/EKF1/EKF1_FHN_STRONG_GLOBAL_ERRORS.csv}\loadedtableIII
        \addplot[gray, only marks, mark=diamond*, mark size=1.7pt] table {\loadedtableIII};
        \addplot[gray, line width=1pt]
        table [y={create col/linear regression={y=errors}}] {\loadedtableIII};
        \xdef\slopeIII{\pgfplotstableregressiona}
        \xdef\bIII{\pgfplotstableregressionb}

        \pgfplotstableread[x=deltas, y=errors, col sep=comma]{numerics/EKF1_2/EKF1_2_FHN_STRONG_GLOBAL_ERRORS.csv}\loadedtableV
        \addplot[coolgrey, only marks, mark=o, mark size=1.7pt] table {\loadedtableV};
        \addplot[coolgrey, line width=1pt]
        table [y={create col/linear regression={y=errors}}] {\loadedtableV};
        \xdef\slopeV{\pgfplotstableregressiona}
        \xdef\bV{\pgfplotstableregressionb}

        \legend{,$\hat{\varepsilon}_{g,eu}\approx \pgfmathparse{exp(\bI)} \xdef\ebI{\pgfmathresult} \pgfmathprintnumber[precision=2]{\ebI} \delta^{\pgfmathprintnumber[precision=2, fixed zerofill]
        {\slopeI}}$,,$\hat{\varepsilon}_{g,\textrm{EKF0}-(2)}\approx \pgfmathparse{exp(\bII)} \xdef\ebII{\pgfmathresult} \pgfmathprintnumber[precision=2]{\ebII} \delta^{\pgfmathprintnumber[precision=2, fixed zerofill]
        {\slopeII}}$,,$\hat{\varepsilon}_{g,\textrm{EKF1}-(2)}\approx \pgfmathparse{exp(\bIII)} \xdef\ebIII{\pgfmathresult} \pgfmathprintnumber[precision=2]{\ebIII}\delta^{\pgfmathprintnumber[precision=2, fixed zerofill]
        {\slopeIII}}$,,$\hat{\varepsilon}_{g,\textrm{EKF1}-(3)}\approx \pgfmathparse{exp(\bV)} \xdef\ebV{\pgfmathresult} \pgfmathprintnumber[precision=2]{\ebV}\delta^{\pgfmathprintnumber[precision=2, fixed zerofill]
        {\slopeV}}$};

    \end{axis}
\end{tikzpicture}
        \end{subfigure}
        \begin{subfigure}[b]{0.47\textwidth}
            \centering
            \begin{tikzpicture}[scale=0.58, font=\LARGE]
    \begin{axis}
        [
        legend cell align={left},
        xlabel={$\delta$},
        ylabel={$L^2$ errors},
        title={Weak error comparisons},
        legend style={
            nodes={scale=0.8, transform shape},
            at={(0,1)},
            anchor=north west,
            fill opacity=0.8, 
            text opacity =1},
        grid=both,
        xmode=log,
        ymode=log,
        xmax=0.05,
        ymin=1e-5, ymax=1,
        ]

        \pgfplotstableread[x=deltas, y=errors, col sep=comma]{numerics/EM/EM_FHN_WEAK_GLOBAL_ERRORS.csv}\loadedtableI
        \addplot[black, only marks, mark=triangle*, mark size=1.7pt] table {\loadedtableI};
        \addplot[black, dashed, line width=1pt]
        table [y={create col/linear regression={y=errors}}] {\loadedtableI};
        \xdef\slopeI{\pgfplotstableregressiona}
        \xdef\bI{\pgfplotstableregressionb}

        \pgfplotstableread[x=deltas, y=errors, col sep=comma]{numerics/EKF0/EKF0_FHN_WEAK_GLOBAL_ERRORS.csv}\loadedtableII
        \addplot[darkgray, only marks, mark=square*, mark size=1.7pt] table {\loadedtableII};
        \addplot[darkgray, dashed, line width=1pt]
        table [y={create col/linear regression={y=errors}}] {\loadedtableII};
        \xdef\slopeII{\pgfplotstableregressiona}
        \xdef\bII{\pgfplotstableregressionb}

        \pgfplotstableread[x=deltas, y=errors, col sep=comma]{numerics/EKF1/EKF1_FHN_WEAK_GLOBAL_ERRORS.csv}\loadedtableIII
        \addplot[gray, only marks, mark=diamond*, mark size=1.7pt] table {\loadedtableIII};
        \addplot[gray, dashed, line width=1pt]
        table [y={create col/linear regression={y=errors}}] {\loadedtableIII};
        \xdef\slopeIII{\pgfplotstableregressiona}
        \xdef\bIII{\pgfplotstableregressionb}

        \pgfplotstableread[x=deltas, y=errors, col sep=comma]{numerics/EKF1_2/EKF1_2_FHN_WEAK_GLOBAL_ERRORS.csv}\loadedtableV
        \addplot[color=coolgrey, only marks, mark=o, mark size=1.7pt] table {\loadedtableV};
        \addplot[coolgrey, dashed, line width=1pt]
        table [y={create col/linear regression={y=errors}}] {\loadedtableV};
        \xdef\slopeV{\pgfplotstableregressiona}
        \xdef\bV{\pgfplotstableregressionb}

        \legend{,$\hat{\varepsilon}_{w,eu}\approx\pgfmathparse{exp(\bI)} \xdef\ebI{\pgfmathresult} \pgfmathprintnumber[precision=2]{\ebI}\delta^{\pgfmathprintnumber[precision=2, fixed zerofill]
        {\slopeI}}$,,$\hat{\varepsilon}_{w,\textrm{EKF0}-(2)}\approx\pgfmathparse{exp(\bII)} \xdef\ebII{\pgfmathresult} \pgfmathprintnumber[precision=2]{\ebII}\delta^{\pgfmathprintnumber[precision=2, fixed zerofill]
        {\slopeII}}$,,$\hat{\varepsilon}_{w,\textrm{EKF1}-(2)}\approx\pgfmathparse{exp(\bIII)} \xdef\ebIII{\pgfmathresult} \pgfmathprintnumber[precision=2]{\ebIII}\delta^{\pgfmathprintnumber[precision=2, fixed zerofill]
        {\slopeIII}}$,,$\hat{\varepsilon}_{w,\textrm{EKF1}-(3)}\approx\pgfmathparse{exp(\bV)} \xdef\ebV{\pgfmathresult} \pgfmathprintnumber[precision=2]{\ebV}\delta^{\pgfmathprintnumber[precision=2, fixed zerofill]
        {\slopeV}}$};

    \end{axis}
\end{tikzpicture}
        \end{subfigure}
        \caption{Fitzhugh--Nagumo model. Left: strong global errors for the Euler-Maruyama scheme with parameter $\delta$ and Algorithms~\ref{alg:ekf0-scheme2} and~\ref{alg:ekf0-scheme3}. Right: similar to the left but with weak global errors.}
        \label{fig:Fig4}
    \end{figure}
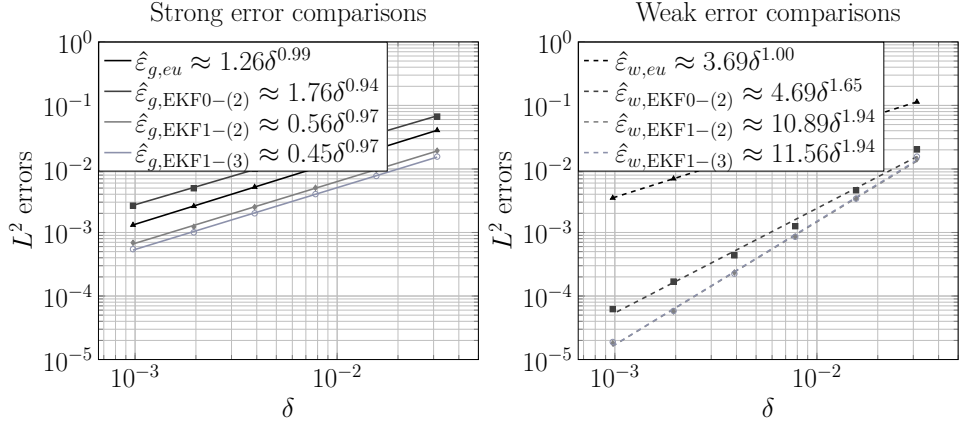
    We draw the following empirical conclusions:
    \begin{enumerate}
    [label=(\roman*)]
        \item The Algorithms~\ref{alg:ekf0-scheme2} and~\ref{alg:ekf0-scheme3} using EKF0 display $1.5$ strong local convergence and $1.0$ strong global convergence, confirming Theorems~\ref{th:ekf0-scheme2} and~\ref{th:ekf0-scheme3}. This is shown by the solid line regressions in the left plots of Figures~\ref{fig:Fig1} and~\ref{fig:Fig2}.
        \item The above conclusion appears to be valid with EKF1 used instead of EKF0, as shown by the right plots in Figures~\ref{fig:Fig1} and~\ref{fig:Fig2}. The convergence orders are the same but EKF1 gives slightly better performances for the multiplicative factor. See the right plot of Figure~\ref{fig:Fig4} for a comparison of the previous schemes as well as the Euler-Maruyama scheme with a step size $\delta$.
        \item All the above schemes, that is EKF0 and EKF1 versions of Algorithms~\ref{alg:ekf0-scheme2} and~\ref{alg:ekf0-scheme3}, exhibit high-order weak convergences, with EKF1 having a higher speed of convergence compared to EKF0. See the dashed line regressions in the right plots of Figures~\ref{fig:Fig1},~\ref{fig:Fig2}, and Figure~\ref{fig:Fig4}.
        \item The marginalised version of the parabola-ODE given by Algorithm~\ref{alg:ekf0-marginal-scheme} with EKF0 linearisation exhibits low-order weak convergence, this can be seen with the grey line regression in the left plot of Figure~\ref{fig:Fig3}. The EKF1 version of Algorithm~\ref{alg:ekf0-marginal-scheme} however does not converge weakly: we observe biased means, see Figure~\ref{fig:Fig3}.
    \end{enumerate}

    \subsection{Uncertainty Calibration}
    In this section, we illustrate the calibration uncertainty methodology proposed in Section~\ref{section:calibration}.
    We compute the maximum-likelihood estimate for the diffusion coefficient $\eta^2$ using~\eqref{eq:mle_estmiate_for_the_diffusion_coefficient} from Proposition~\ref{prop:ml_post_hoc}.
    We first compute a pathwise solution using Algorithm~\ref{alg:ekf0-scheme3} with $\eta^2=1.0$, which is then used to compute the post hoc estimate $\hat{\eta}^2$~\eqref{eq:mle_estmiate_for_the_diffusion_coefficient}. The calibrated solution is then obtained by setting $\eta^2=\hat{\eta}^2$ in Algorithm~\ref{alg:ekf0-scheme3}.
    See Figure~\ref{fig:FigCalibration} for the comparison between the EKF0/1 approximations given by Algorithm~\ref{alg:ekf0-scheme3} with $\eta = \sqrt{10}$ and $\eta=\hat{\eta}_{\textup{MLE}}$.
    The step size is $\delta=0.02$, the fine solution (FS) denoted by $(X_{t_k}^{\textup{em}})_{0\leq k\leq K}$ is obtained with the EM scheme with step size $\delta^2$.
    In order to help with visualisation, in addition, we report the marginal Z-scores associated with both components in Figure~\ref{fig:FigCalibrationZscore}, computed as follows:
    \begin{equation}\label{eq:z-scores}
        \frac{X_{t_k}^{em, (i)}-m_{k-1}^{(i)}(t_k)}{\sqrt{P_{k-1}^{(i,i)}(t_k)}},
    \end{equation}
    for component $i\in\{1, 2\}$ and time step $t_k$.
    \begin{figure}[t]
        \centering
        \begin{subfigure}[b]{0.49\textwidth}
            \begin{tikzpicture}[scale=0.58, font=\LARGE]
    \begin{axis}
        [
        legend cell align={left},
        xlabel={$t$},
        ylabel={},
        title={EKF0},
        legend style={
            nodes={scale=0.8, transform shape},
            at={(0,1)},
            anchor=north west,
            fill opacity=0.8, 
            text opacity =1},
        grid=both,
        xmin=0, xmax=1,
        ymin=-1.5, ymax=0.7,
        error bars/y dir=both, %
        error bars/y explicit  %
        ]

        \pgfplotstableread[x=t, y=mean, col sep=comma]{numerics_review/UC_EKF0_Alg3_FHN/EM_FHN_calibrated_fine_0_50_50.csv}\loadedtableI
        \addplot[black, line width=1pt] table {\loadedtableI};
        
        \pgfplotstableread[x=t, y=mean, col sep=comma]{numerics_review/UC_EKF0_Alg3_FHN/EM_FHN_calibrated_fine_1_50_50.csv}\loadedtableII
        \addplot[darkgray, line width=1pt] table {\loadedtableII};
        
        \pgfplotstableread[x=t, y=mean, y error=error, col sep=comma]{numerics_review/UC_EKF0_Alg3_FHN/EKF0_2_FHN_uncalibrated_UC_0_50_50.csv}\loadedtableIII
        \addplot[gray, line width=1pt, error bars/.cd] table [y error=error] {\loadedtableIII};
        
        \pgfplotstableread[x=t, y=mean, y error=error, col sep=comma]{numerics_review/UC_EKF0_Alg3_FHN/EKF0_2_FHN_uncalibrated_UC_1_50_50.csv}\loadedtableIV
        \addplot+[coolgrey, line width=1pt, ] table [y error=error] {\loadedtableIV};

        \legend{$1^{\textup{st}}$ coord. FS, $2^{\textup{nd}}$ coord. FS, $1^{\textup{st}}$ coord., $2^{\textup{nd}}$ coord.};
        
    \end{axis}
\end{tikzpicture}
        \end{subfigure}
        \begin{subfigure}[b]{0.49\textwidth}
            \begin{tikzpicture}[scale=0.58, font=\LARGE]
    \begin{axis}
        [
        legend cell align={left},
        xlabel={$t$},
        ylabel={},
        title={EKF1},
        legend style={
            nodes={scale=0.8, transform shape},
            at={(0,1)},
            anchor=north west,
            fill opacity=0.8, 
            text opacity =1},
        grid=both,
        xmin=0, xmax=1,
        ymin=-1.5, ymax=0.7,
        error bars/y dir=both, %
        error bars/y explicit  %
        ]

        \pgfplotstableread[x=t, y=mean, col sep=comma]{numerics_review/UC_EKF1_Alg3_FHN/EM_FHN_calibrated_fine_0_50_50.csv}\loadedtableI
        \addplot[black, line width=1pt] table {\loadedtableI};
        
        \pgfplotstableread[x=t, y=mean, col sep=comma]{numerics_review/UC_EKF1_Alg3_FHN/EM_FHN_calibrated_fine_1_50_50.csv}\loadedtableII
        \addplot[darkgray, line width=1pt] table {\loadedtableII};
        
        \pgfplotstableread[x=t, y=mean, y error=error, col sep=comma]{numerics_review/UC_EKF1_Alg3_FHN/EKF1_2_FHN_uncalibrated_UC_0_50_50.csv}\loadedtableIII
        \addplot[gray, line width=1pt, error bars/.cd] table [y error=error] {\loadedtableIII};
        
        \pgfplotstableread[x=t, y=mean, y error=error, col sep=comma]{numerics_review/UC_EKF1_Alg3_FHN/EKF1_2_FHN_uncalibrated_UC_1_50_50.csv}\loadedtableIV
        \addplot[coolgrey, line width=1pt, ] table [y error=error] {\loadedtableIV};

        \legend{$1^{\textup{st}}$ coord. FS, $2^{\textup{nd}}$ coord. FS, $1^{\textup{st}}$ coord., $2^{\textup{nd}}$ coord.};
        
    \end{axis}
\end{tikzpicture}
        \end{subfigure}
        \medskip
        \begin{subfigure}[b]{0.49\textwidth}
            \begin{tikzpicture}[scale=0.58, font=\LARGE]
    \begin{axis}
        [
        legend cell align={left},
        xlabel={$t$},
        ylabel={},
        title={EKF0},
        legend style={
            nodes={scale=0.8, transform shape},
            at={(0,1)},
            anchor=north west,
            fill opacity=0.8, 
            text opacity =1},
        grid=both,
        xmin=0, xmax=1,
        ymin=-1.5, ymax=0.7,
        error bars/y dir=both, %
        error bars/y explicit  %
        ]

        \pgfplotstableread[x=t, y=mean, col sep=comma]{numerics_review/UC_EKF0_Alg3_FHN/EM_FHN_calibrated_fine_0_50_50.csv}\loadedtableI
        \addplot[black, line width=1pt] table {\loadedtableI};
        
        \pgfplotstableread[x=t, y=mean, col sep=comma]{numerics_review/UC_EKF0_Alg3_FHN/EM_FHN_calibrated_fine_1_50_50.csv}\loadedtableII
        \addplot[darkgray, line width=1pt] table {\loadedtableII};
        
        \pgfplotstableread[x=t, y=mean, y error=error, col sep=comma]{numerics_review/UC_EKF0_Alg3_FHN/EKF0_2_FHN_calibrated_UC_0_50_50.csv}\loadedtableIII
        \addplot[gray, line width=1pt, error bars/.cd] table [y error=error] {\loadedtableIII};
        
        \pgfplotstableread[x=t, y=mean, y error=error, col sep=comma]{numerics_review/UC_EKF0_Alg3_FHN/EKF0_2_FHN_calibrated_UC_1_50_50.csv}\loadedtableIV
        \addplot+[coolgrey, line width=1pt, ] table [y error=error] {\loadedtableIV};

        \legend{$1^{\textup{st}}$ coord. FS, $2^{\textup{nd}}$ coord. FS, $1^{\textup{st}}$ coord., $2^{\textup{nd}}$ coord.};
        
    \end{axis}
\end{tikzpicture}
        \end{subfigure}
        \begin{subfigure}[b]{0.49\textwidth}
            \begin{tikzpicture}[scale=0.58, font=\LARGE]
    \begin{axis}
        [
        legend cell align={left},
        xlabel={$t$},
        ylabel={},
        title={EKF1},
        legend style={
            nodes={scale=0.8, transform shape},
            at={(0,1)},
            anchor=north west,
            fill opacity=0.8, 
            text opacity =1},
        grid=both,
        xmin=0, xmax=1,
        ymin=-1.5, ymax=0.7,
        error bars/y dir=both, %
        error bars/y explicit  %
        ]

        \pgfplotstableread[x=t, y=mean, col sep=comma]{numerics_review/UC_EKF1_Alg3_FHN/EM_FHN_calibrated_fine_0_50_50.csv}\loadedtableI
        \addplot[black, line width=1pt] table {\loadedtableI};
        
        \pgfplotstableread[x=t, y=mean, col sep=comma]{numerics_review/UC_EKF1_Alg3_FHN/EM_FHN_calibrated_fine_1_50_50.csv}\loadedtableII
        \addplot[darkgray, line width=1pt] table {\loadedtableII};
        
        \pgfplotstableread[x=t, y=mean, y error=error, col sep=comma]{numerics_review/UC_EKF1_Alg3_FHN/EKF1_2_FHN_calibrated_UC_0_50_50.csv}\loadedtableIII
        \addplot[gray, line width=1pt, error bars/.cd] table [y error=error] {\loadedtableIII};
        
        \pgfplotstableread[x=t, y=mean, y error=error, col sep=comma]{numerics_review/UC_EKF1_Alg3_FHN/EKF1_2_FHN_calibrated_UC_1_50_50.csv}\loadedtableIV
        \addplot[coolgrey, line width=1pt, ] table [y error=error] {\loadedtableIV};

        \legend{$1^{\textup{st}}$ coord. FS, $2^{\textup{nd}}$ coord. FS, $1^{\textup{st}}$ coord., $2^{\textup{nd}}$ coord.};
        
    \end{axis}
\end{tikzpicture}
        \end{subfigure}
        \caption{Algorithm~\ref{alg:ekf0-scheme3}, Fitzhugh--Nagumo model. Top, prior with $\eta=\sqrt{10}$. Bottom, calibrated prior with $\eta=\hat{\eta}_{\textup{MLE}}$~\eqref{eq:mle_estmiate_for_the_diffusion_coefficient}. For the two coordinates, the fine solution (FS) sample path obtained via EM and the Gaussian approximation is given by its mean and its standard deviation.}
        \label{fig:FigCalibration}
    \end{figure}
    \begin{figure}[t]
        \centering
        \begin{subfigure}[b]{0.49\textwidth}
            \begin{tikzpicture}[scale=0.58, font=\LARGE]
    \begin{axis}
        [
        legend cell align={left},
        xlabel={$t$},
        ylabel={},
        title={EKF0},
        legend style={
            nodes={scale=0.8, transform shape},
            at={(0,1)},
            anchor=north west,
            fill opacity=0.8, 
            text opacity =1},
        grid=both,
        xmin=0, xmax=1,
        ymin=-9, ymax=2,
        error bars/y dir=both, %
        error bars/y explicit  %
        ]

        \pgfplotstableread[x=t, y=mean, col sep=comma]{numerics_review/UC_EKF0_Alg3_FHN/EKF0_2_FHN_uncalibrated_zscore_0_50_50.csv}\loadedtableI
        \addplot[black, line width=1pt, dashed] table {\loadedtableI};

        \pgfplotstableread[x=t, y=mean, col sep=comma]{numerics_review/UC_EKF0_Alg3_FHN/EKF0_2_FHN_uncalibrated_zscore_1_50_50.csv}\loadedtableII
        \addplot[gray, line width=1pt, dashed] table {\loadedtableII};

        \pgfplotstableread[x=t, y=mean, col sep=comma]{numerics_review/UC_EKF0_Alg3_FHN/EKF0_2_FHN_calibrated_zscore_0_50_50.csv}\loadedtableIII
        \addplot[black, line width=1pt] table {\loadedtableIII};

        \pgfplotstableread[x=t, y=mean, col sep=comma]{numerics_review/UC_EKF0_Alg3_FHN/EKF0_2_FHN_calibrated_zscore_1_50_50.csv}\loadedtableIV
        \addplot[gray, line width=1pt] table {\loadedtableIV};
        
        \legend{$1^{\textup{st}}$ coord. $\eta=\sqrt{10}$, $2^{\textup{nd}}$ coord. $\eta=\sqrt{10}$, $1^{\textup{st}}$ coord. $\eta=\hat{\eta}_{\textup{MLE}}$, $2^{\textup{nd}}$ coord. $\eta=\hat{\eta}_{\textup{MLE}}$};
        
    \end{axis}
\end{tikzpicture}
        \end{subfigure}
        \begin{subfigure}[b]{0.49\textwidth}
            \begin{tikzpicture}[scale=0.58, font=\LARGE]
    \begin{axis}
        [
        legend cell align={left},
        xlabel={$t$},
        ylabel={},
        title={EKF1},
        legend style={
            nodes={scale=0.8, transform shape},
            at={(0,1)},
            anchor=north west,
            fill opacity=0.8, 
            text opacity =1},
        grid=both,
        xmin=0, xmax=1,
        ymin=-2, ymax=2,
        error bars/y dir=both, %
        error bars/y explicit  %
        ]

        \pgfplotstableread[x=t, y=mean, col sep=comma]{numerics_review/UC_EKF1_Alg3_FHN/EKF1_2_FHN_uncalibrated_zscore_0_50_50.csv}\loadedtableI
        \addplot[black, line width=1pt, dashed] table {\loadedtableI};

        \pgfplotstableread[x=t, y=mean, col sep=comma]{numerics_review/UC_EKF1_Alg3_FHN/EKF1_2_FHN_uncalibrated_zscore_1_50_50.csv}\loadedtableII
        \addplot[gray, line width=1pt, dashed] table {\loadedtableII};

        \pgfplotstableread[x=t, y=mean, col sep=comma]{numerics_review/UC_EKF1_Alg3_FHN/EKF1_2_FHN_calibrated_zscore_0_50_50.csv}\loadedtableIII
        \addplot[black, line width=1pt] table {\loadedtableIII};

        \pgfplotstableread[x=t, y=mean, col sep=comma]{numerics_review/UC_EKF1_Alg3_FHN/EKF1_2_FHN_calibrated_zscore_1_50_50.csv}\loadedtableIV
        \addplot[gray, line width=1pt] table {\loadedtableIV};
        
        \legend{$1^{\textup{st}}$ coord. $\eta=\sqrt{10}$, $2^{\textup{nd}}$ coord. $\eta=\sqrt{10}$, $1^{\textup{st}}$ coord. $\eta=\hat{\eta}_{\textup{MLE}}$, $2^{\textup{nd}}$ coord. $\eta=\hat{\eta}_{\textup{MLE}}$};
        
    \end{axis}
\end{tikzpicture}
        \end{subfigure}
        \caption{Algorithm~\ref{alg:ekf0-scheme3}, Fitzhugh--Nagumo model. Dashed lines, prior with $\eta=\sqrt{10}$. Solid lines, calibrated prior with $\eta=\hat{\eta}_{\textup{MLE}}$~\eqref{eq:mle_estmiate_for_the_diffusion_coefficient}. For the two coordinates, $Z$-scores of the fine solution associated with the posterior Gaussian approximation.}
        \label{fig:FigCalibrationZscore}
    \end{figure}
    It can be seen that the non-calibrated posterior solutions (top row) produce errors that oscillate outside their reported uncertainties while the calibrated ones largely do not, hinting to a better coverage of the posterior uncertainty.

    \section{Discussion}
    \label{section:discussion}
    In this article, we have presented a probabilistic numerics approach to characterising the discretisation error in additive-noise SDEs. We believe that this paves the way for the development of further probabilistic numerics techniques for sampling from SDE solutions.
    At its core, our methods rely on computing a corrected posterior distribution over the sample path solution of the SDE.
    We use piecewise differentiable approximations of the Brownian motion on top of which we apply existing Gaussian ODE filtering methods~\citep{HennigKerstingTronarpSarkka2019, HennigTronarpSarkka2021}, allowing for tractable estimation of approximated Gaussian transition densities.
    To quantify uncertainty, our approach increases computational complexity compared to traditional pointwise SDE simulators, typically costing $d^3$ due to matrix inversions appearing in the Kalman-filtering routine and potentially (for the EKF1 variant of our method here) requiring gradient evaluations.
    However, and similarly to Bayesian ODE filters, it can benefit from solutions built either on independence assumptions or Kronecker structures on the model~\citep{kramer2022Prob}, thus making the inference computationally almost as efficient as the EM scheme or other standard SDE samplers.

    To summarise, we have introduced three key methods.
    The first one, the Gaussian SDE filter, consists of drawing a local Brownian approximation, applying a single Kalman filter step for each ODE, and then sampling from the posterior solution.
    The second one, the Mixture Gaussian SDE filter is a variation of the previous one. Instead of sampling at each time step a solution from the posterior, we propagate the posterior uncertainties through the Kalman filter, leading to a mixture of Gaussian posteriors given a sample Brownian path.
    The third one, the Marginalised Gaussian SDE filter, leverages the parabola approximation of the Brownian motion introduced in~\citep{Foster2020OptimalPolyBrownian} and incorporates the Gaussian coefficients as part of the observed variables, thus marginalising the Brownian path along the SDE solution and resulting in a joint transition density over the solution and the Brownian path.

    We have proven empirically and theoretically that our methods retain good pathwise and moment properties.
    The first two schemes in particular are provably 1.0 strongly convergent in the specific instance where EKF0 is used on top of the parabola approximation and for an IOUP prior. Moreover, numerical experiments highlight a weak convergence order $> 1.5$ that is yet to be analysed in future work.
    We anticipate that a similar proof can be adapted for EKF1 (Algorithm~\ref{alg:ekf0-scheme2}), which empirically exhibits similar order convergence as EKF0 (respectively strong local order $1.5$, strong global order $1$, and weak order $2$).
    Moreover, we believe Assumption~\ref{assumption:reg_of_sde_flow} of bounded drift functions can be weakened: indeed the drift function in the FHN model is not bounded, does not admit bounded derivative, and is not Lipschitz; this did not prevent the method from working in this case.
    The third, marginalised, scheme appears numerically to be weakly convergent with low-order convergence when EKF0 is used. However, when EKF1 is used to propagate uncertainty on the solution, the resulting trajectory estimate suffers from non-negligible biased moment estimations, likely resulting from an error accumulation that is absent from EKF0 due to its more radical forgetful initialisation.

    Perhaps the main drawback of our method lies in its \emph{underconfident} uncertainty estimates, which is a problem inherited from probabilistic numerics for ordinary differential equations, specifically when using EKF1-type linearisations~\citet{Bosch2020CalibratedAP}.
    The question of better calibrating uncertainty is an active topic of research in the literature on Gaussian-approximated ordinary differential equations~\citep[see, e.g.]{Kersting2016active,Bosch2020CalibratedAP}, and improvements made in this direction are likely to translate immediately to our setting provided the underlying prior/inference methods are compatible with a piecewise approximation approach.

    While our current algorithm leverages the specific piecewise Brownian approximation of~\citet{Foster2020OptimalPolyBrownian}, it is worth noting that the flexibility of our approach allows for exploring the same framework with different approximations, including higher-order ones, which could lead to performance improvements. Similarly, a large number of extensions of the Gaussian ODE framework have been proposed in the recent literature~\citep[see, e.g.,][]{kramer2022Prob,bosch2023probabilistic}, which could be used as a drop-in replacement for EKF0 or EKF1.
    Furthermore, our Bayesian approach allows for potential improvements of the sampler through the use of adapted solution priors, and uncertainty calibration methods.
    For example, using periodic priors for periodic SDEs, or Matérn priors for differentiable SDE paths, could enhance the numerical performances of the samplers or, similar to~\citet{HennigKerstingTronarpSarkka2019}, calibrating the prior noise could also lead the SDE sampler to report appropriate uncertainty. Additionally, better Gaussian integrators, such as iterative and implicit methods~\citep{HennigTronarpSarkka2021}, are likely directly applicable to our work. We leave the study of their empirical and theoretical properties to future work.

    \section*{Acknowledgment}
    The authors would like to thank two anonymous referees as well as an Associate Editor for their constructive feedback which prompted the discussion on uncertainty calibration and improved the quality of the manuscript overall.

    \section*{Individual contributions}
    The original idea for this article is due to AC. The methodology was primarily developed by both AC and YLF with inputs from SS. The implementation and experimental evaluation are due to YLF with inputs from AC. The theory is due to YLF building on ideas from AC. The original version of this article was written by YLF with feedback and corrections from AC, after which all authors reviewed the manuscript. The project was coordinated by SS.

        \section{Supplementary Material}
    The supplementary material contains proofs of our main results (Theorem~\ref{th:ekf0-scheme2} and~\ref{th:ekf0-scheme3}) in Sections~\ref{app:ekf0-scheme2} and~\ref{app:ekf0-scheme3}, respectively, of auxiliary lemmas necessary for these proofs in Section~\ref{app:aux}, and further details regarding the multivariate formulation for the marginalisation procedure of Section~\ref{section:extension2} in Section~\ref{appendix:multivariate_case_marginal_brownian}.

    For the sake of notation simplicity, we omit the superscript dependency on $\delta$ of the Brownian approximation $\beta^{\delta}$, which we simply denote by $\beta$. We denote by $\dot{\beta}$ the time-derivative of $\beta$ and similarly for $\ddot{\beta}$. For some real-valued functions $f$ and $g$, we write $f(\delta)=O(g(\delta))$ if $\abs{f(\delta)/g(\delta)} \leq M$ for some constant $M$ as $\delta$ goes to $0$. This is the standard big-$O$ notation. We write $f(\delta) = O_{\calL^2}(g(\delta))$ if $\normLtwo{f(\delta)}=O(g(\delta))$.

    \label{appendix:ekf0_proof}
    We consider the univariate case $d=1$, \textit{mutatis mutandis}, and the proofs are still valid in the multidimensional case.

    \subsection{Auxiliary results}\label{app:aux}
    In this section, we prove intermediary results that will be useful in proving Theorems~\ref{th:ekf0-scheme2} and~\ref{th:ekf0-scheme3}. Lemma~\ref{lemma:diff_integral} is concerned with an $\calL^2$-expansion of the term $\int_{0}^{\delta}\sigma(u)\dd \beta(u)$ which is involved in the local steps given by Algorithms~\ref{alg:ekf0-scheme2},~\ref{alg:ekf0-scheme3}, and the local expansion of the parabola-ODE solution. In Lemma~\ref{lemma:kalman_gain}, we compute intermediary quantities needed to compute the EKF0 step. In Lemma~\ref{lemma:localekf0}, we derive the local step of Algorithm~\ref{alg:ekf0-scheme2}, in Lemma~\ref{lemma:local_step_parabola_ode}, we compute a local expansion of the exact parabola-ODE solution.

    \begin{lemma}
        \label{lemma:diff_integral}
        Let $\sigma$ be a continuously differentiable function over $[0, T]$.
        Let $\beta$ be the parabola approximation of a Brownian motion over $[0, \delta]$.
        Let $A_{-} = \dot{\beta}(0)$ and $A_{+} = \dot{\beta}(\delta)$.
        Then for any $0 < \delta \leq T$
        \begin{equation*}
            \biggnormLtwo{\int_{0}^{\delta}\sigma(u)\dd \beta(u)-\frac{\delta}{2}(\sigma(0)A_{-}+\sigma(\delta)A_{+})}=O(\delta^{1.5}),
        \end{equation*}
        which can we written as
        \begin{equation*}
            \int_{0}^{\delta}\sigma(u)\dd \beta(u) = \frac{\delta}{2}(\sigma(0)A_{-}+\sigma(\delta)A_{+}) + O_{\calL^2}(\delta^{1.5}).
        \end{equation*}
        \begin{proof}
            Let $\delta > 0$.
            We denote by $\beta$ the parabola approximation over the interval $[0,\delta]$, that is, for $u\in [0, \delta]$,
            \begin{equation*}
                \beta(u) = \frac{B_{\delta} u }{\delta} + \frac{\sqrt{6} I_{\delta} u }{\delta}\left(\frac{u}{\delta}-1\right),
            \end{equation*}
            where $B_{\delta}$ and $I_{\delta}$ are independent central normal variables with variance $\delta$ and $\delta/2$, respectively.
            Let $J = \int_{0}^{\delta}\sigma(u)\dd \beta(u)$.
            Since $\dot{\beta}$ is a linear function, we have
            \begin{equation}
                \begin{split}
                    \label{eq:proof0}
                    \dot{\beta}(u) &= \dot{\beta}(0) + u\ddot{\beta}(0) = A_{-} + u \ddot{\beta},\\
                    &= \dot{\beta}(\delta) + (u-\delta)\ddot{\beta}(\delta) = A_{+} + (u-\delta)\ddot{\beta},
                \end{split}
            \end{equation}
            where
            \begin{equation}
                \label{eq:def_A_bpp}
                A_{\pm} = \frac{B_{\delta}}{\delta} \pm \frac{\sqrt{6}I_{\delta}}{\delta}, \quad \ddot{\beta} = \frac{2\sqrt{6}I_{\delta}}{\delta^2}.
            \end{equation}
            Now, let us tackle the diffusion part in $J$.
            For any $u\in [0,\delta/2]$, the mean-value theorem ensures the existence of $\theta_{u}^{1}\in [0,u]$ such that
            \begin{equation*}
                \sigma(u) = \sigma(0) + u\sigma'(\theta_{u}^{1}).
            \end{equation*}
            Similarly for any $u\in [\delta/2, \delta]$, there exists $\theta_{u}^{2}\in [u, \delta]$ such that
            \begin{equation*}
                \sigma(u) = \sigma(\delta) + (u-\delta)\sigma'(\theta_{u}^{2}).
            \end{equation*}
            Using~\eqref{eq:proof0} and the previous applications of the mean-value theorem, we have
            \begin{equation}
                \label{eq:equality_integral_diff}
                \begin{split}
                    J &= \int_{0}^{\delta/2}(\sigma(0)+u\sigma'(\theta_{u}^{1}))(A_{-}+u\ddot{\beta})\dd u + \int_{\delta/2}^{\delta}(\sigma(\delta)+(u-\delta)\sigma'(\theta_{u}^{2}))(A_{+} + (u-\delta)\ddot{\beta})\dd u\\
                    &= \frac{\delta}{2}(\sigma(0)A_{-} + \sigma(\delta)A_{+}) + R_1 + R_2 + R_3,
                \end{split}
            \end{equation}
            where $R_1$, $R_2$, and $R_3$ are defined as
            \begin{equation}
                \label{eq:def_rests}
                \begin{split}
                    R_1 &= \ddot{\beta} \cdot \left[\int_{0}^{\delta/2} u^2\sigma'(\theta_{u}^{1})\dd u + \int_{\delta/2}^{\delta} (u-\delta)^2\sigma'(\theta_{u}^{2})\dd u\right],\\
                    R_2 &= A_{-}\int_{0}^{\delta/2} u \sigma'(\theta_{u}^{1}) \dd u + A_{+}\int_{\delta/2}^{\delta} (u-\delta)\sigma'(\theta_{u}^{2})\dd u,\\
                    R_3 &= \ddot{\beta} \cdot \left[\int_{0}^{\delta/2}\sigma(0)u\dd u + \int_{\delta/2}^{\delta}\sigma(\delta)(u-\delta)\dd u\right].
                \end{split}
            \end{equation}
            Using~\eqref{eq:def_A_bpp}, we have
            \begin{equation}
                \label{eq:bound_l2_norm}
                \begin{split}
                    \normLtwo{A_{\pm}} = \frac{2}{\sqrt{\delta}} = O( \delta^{-0.5}),\,\,
                    \normLtwo{\ddot{\beta}} = \frac{2\sqrt{3}}{\delta^{1.5}} = O(\delta^{-1.5}).
                \end{split}
            \end{equation}
            Since $\sigma$ is continuously differentiable, its derivative can be uniformly bounded on the closed interval $[0, T]$,~\eqref{eq:def_rests} combined to~\eqref{eq:bound_l2_norm} leads to for the first two terms
            \begin{equation}
                \label{eq:bound_l2_norm2}
                \begin{split}
                    \normLtwo{R_1} &\leq \max_{t\in[0, T]}\abs{\sigma'(t)} \cdot \normLtwo{\ddot{\beta}}\cdot\left(\int_{0}^{\delta/2}u^2\dd u+\int_{\delta/2}^{\delta}(u-\delta)^2\dd u\right)\\
                    &=\max_{t\in[0, T]}\abs{\sigma'(t)} \cdot \normLtwo{\ddot{\beta}}\cdot\delta^3/12= O(\delta^{1.5}),\\
                    \normLtwo{R_2} &\leq \max_{t\in[0,T]}\abs{\sigma'(t)} \cdot \left(\normLtwo{A_{-}}\int_{0}^{\delta/2} u\dd u+\normLtwo{A_{+}}\int_{\delta/2}^{\delta}(\delta-u)\dd u\right) \\
                    &= \max_{t\in[0,T]}\abs{\sigma'(t)} \cdot (\normLtwo{A_{-}}+\normLtwo{A_{+}}) \cdot \delta^2/8 = O(\delta^{1.5}).
                \end{split}
            \end{equation}
            For the last term $R_3$, using the mean-value theorem, there exists $\theta^3\in [0,\delta]$ such that $\sigma(\delta) = \sigma(0)+\delta\sigma'(\theta^3)$.
            Thus replacing $\sigma(\delta)$ in $R_3$, we have
            \begin{equation*}
                \begin{split}
                    R_3 &= \ddot{\beta}\cdot\left[\int_{0}^{\delta/2}\sigma(0)u\dd u + \int_{\delta/2}^{\delta}\sigma(\delta)(u-\delta)\dd u\right]\\
                    &= \ddot{\beta}\cdot\left[\int_{0}^{\delta/2}\sigma(0)u\dd u + \int_{\delta/2}^{\delta}\sigma(0) (u-\delta) \dd u + \int_{\delta/2}^{\delta}\delta \sigma'(\theta^3)(u-\delta)\dd u\right] \\
                    &= \ddot{\beta}\cdot(0 - \sigma'(\theta^3)\delta^3/8).
                \end{split}
            \end{equation*}
            Using~\eqref{eq:bound_l2_norm}, we have
            \begin{equation}
                \label{eq:bound_l2_norm3}
                \normLtwo{R_3} \leq \max_{t\in[0, T]}\abs{\sigma'(t)} \cdot \normLtwo{\ddot{\beta}}\cdot\delta^3/8 = O(\delta^{1.5}).
            \end{equation}
            Combining \eqref{eq:equality_integral_diff} with the $\calL^2$ bounds derived above~\eqref{eq:bound_l2_norm2} and~\eqref{eq:bound_l2_norm3}, the lemma follows.
        \end{proof}
    \end{lemma}
    \begin{lemma}[Kalman gain]
        \label{lemma:kalman_gain}
        Assume the initial variance $P(0)$ and the noise $R$ are $0$.
        Then the Kalman-filter gain $K(\delta)$ defined by~\eqref{eq:kalman_filter_update_eq} is equal to
        \begin{equation*}
            K(\delta) = \begin{pmatrix}
                            \delta/2+O(\delta^3) \\
                            1
            \end{pmatrix}.
        \end{equation*}
        \begin{proof}
            Let $\theta \geq 0$, $\eta >0$, and $\delta > 0 $.
            Recall that the integrated Ornstein-Uhlenbeck process (IOUP) with parameter $\theta$ and diffusion $\eta$, has transition $Y(t+\delta)\mid Y(t) \sim \calN(AY(t), Q)$, where $A$ and $Q$ are defined by
            \begin{equation*}
                \begin{split}
                    A = e^{F\delta},\quad Q = \int_{0}^{\delta}e^{Fs}L L^{\top}e^{F^{\top}s}\dd s\label{eq:def_Q},
                \end{split}
            \end{equation*}
            and where $F$ and $L$ are defined as
            \begin{equation*}
                \begin{split}
                    F =
                    \begin{pmatrix}
                        0 & 1       \\
                        0 & -\theta
                    \end{pmatrix},
                    \quad
                    L  =
                    \begin{pmatrix}
                        0 \\
                        \eta
                    \end{pmatrix}.
                \end{split}
            \end{equation*}
            Immediate calculations show that for any $s > 0$,
            \begin{equation}
                \label{eq:exp_matrix}
                e^{Fs} = \begin{pmatrix}
                             1 & (1-e^{-\theta s})/\theta \\
                             0 & e^{-\theta s}
                \end{pmatrix},
            \end{equation}
            where the case $\theta = 0$ is taken as the limit $\theta \to 0$.
            Using~\eqref{eq:def_Q} and~\eqref{eq:exp_matrix}, we have
            \begin{equation}
                \label{eq:dev_Q}
                \begin{split}
                    Q  &= \eta^2 \int_{0}^{\delta}
                    \begin{pmatrix}
                        ((1-e^{-\theta s})/\theta)^2          & e^{-\theta s}(1-e^{-\theta s})/\theta \\
                        e^{-\theta s}(1-e^{-\theta s})/\theta & e^{-2\theta s}
                    \end{pmatrix}\dd s \\
                    &=\eta^2
                    \begin{pmatrix}
                        \delta^3/3-\theta \delta^4/4+O(\delta^5)   & \delta^2/2-\theta \delta^3/2 + O(\delta^4) \\
                        \delta^2/2-\theta \delta^3/2 + O(\delta^4) & \delta - \theta \delta^2 +O(\delta^3)
                    \end{pmatrix}.
                \end{split}
            \end{equation}
            Direct computations using the update equations~\eqref{eq:kalman_filter_update_eq} with \begin{equation}
                                                                                                       P^{-}=Q, R=0, \tilde{H} = H_1 = (0, 1),
            \end{equation}
            and the previous development in $\delta$~\eqref{eq:dev_Q}, give
            \begin{equation*}
                \begin{aligned}[t]
                    S &= H_1P^{-}H_1^{\top}+R\\
                    &=\eta^2[\delta-\theta \delta^2 + O(\delta^3)],
                \end{aligned}
                \qquad
                \begin{aligned}[t]
                    K &= P^{-}H_1^{\top}S^{-1} \\
                    &=\begin{pmatrix}
                          \delta/2+O(\delta^3) \\
                          1
                    \end{pmatrix}.
                \end{aligned}
            \end{equation*}
        \end{proof}
        \begin{remark}
        \label{rem:kalman_gain}
            The previous computations can be generalized to any noise $R$. In particular, let $R=c\delta^2$ for some constant $c>0$, then
            \begin{equation*}
                K = \begin{pmatrix}
                    \delta/2+O(\delta^2)\\
                    1+O(\delta)
                \end{pmatrix}.
            \end{equation*}
        \end{remark}
    \end{lemma}
    \begin{lemma}[Local EKF0 step]
        \label{lemma:localekf0}
        Let $\delta > 0$, and $Y(\delta)\sim\calN(m(\delta),P(\delta))$ be the distribution of the solution after one step of Algorithm~\ref{alg:ekf0-scheme2} with EKF0 linearisation and with initialisation $\wtX_0 = x_0$ and let $\wtX_{\delta}\sim \calL(Y^{(0)}(\delta))$ be drawn from the output distribution.
        Then under Assumption~\ref{assumption:reg_of_sde_flow},
        \begin{equation*}
            \wtX_{\delta} =  x_0 + \delta \mu(x_0, 0) + \frac{\delta}{2}(\sigma(0)A_{-} + \sigma(\delta)A_{+}) +  \tilde{\varepsilon}_1 + \tilde{\nu}_1,
        \end{equation*}
        where $\tilde{\varepsilon}_1$ satisfies
        \begin{equation*}
            \begin{split}
                \normLtwo{\tilde{\varepsilon}_1} = O(\delta^{1.5}),\quad \bbE[\tilde{\varepsilon}_1] = O(\delta^2),
            \end{split}
        \end{equation*}
        and where the constants in front of the big-$O$ terms depends only on uniform bounds on $\mu$ and $\sigma$, and $\tilde{\nu}_1\sim \calN(0, \eta^2\delta^3/12+O(\delta^5))$ is independent from $\tilde{\varepsilon}_1$.
        \begin{proof}
            The initial mean is set to
            \begin{equation*}
                m(0) =
                \begin{pmatrix}
                    x_0 \\
                    \mu(x_0, 0)+\sigma(0)A_{-}
                \end{pmatrix}.
            \end{equation*}
            Using the predictive equations~\eqref{eq:pred_eq_kalman} and~\eqref{eq:exp_matrix}, we have
            \begin{equation}
                \label{eq:predicted_mean}
                \begin{split}
                    m^{-}(\delta) = Am(0) =
                    \begin{pmatrix}
                        x_0 + (\delta +O(\delta^2))\times(\mu(x_0, 0) + \sigma(0)A_{-}) \\
                        (1+O(\delta))\times(\mu(x_0, 0) + \sigma(0)A_{-})
                    \end{pmatrix}.\\
                \end{split}
            \end{equation}
            The update equations~\eqref{eq:kalman_filter_update_eq} give
            \begin{equation*}
                m^{(0)}(\delta) = m^{-(0)}(\delta)-K^{(0)}\times(m^{-(1)}(\delta)-(\mu(m^{-(0)}, \delta)+\sigma(\delta)A_{+})),
            \end{equation*}
            where $K^{(0)}$ is the first component of the Kalman gain, $m^{-(0)}$ and $m^{-(1)}$ are the predicted mean values for the trajectory and the vector field evaluated at the trajectory.
            Using Lemma~\ref{lemma:kalman_gain} on the Kalman gain and replacing the predicted mean values by~\eqref{eq:predicted_mean}, we obtain
            \begin{equation}
                \label{eq:local_step_ekf0}
                \begin{split}
                    m^{(0)}(\delta) & = x_0 + [\delta+O(\delta^2)]\times[\mu(x_0, 0) + \sigma(0)A_{-}] - [\delta/2+O(\delta^3)]\\
                    &\quad \times \Bigg((1+O(\delta))\times(\mu(x_0, 0) + \sigma(0)A_{-})\\
                    &\qquad- \mu\bigg(x_0 + [\delta +O(\delta^2)]\times[\mu(x_0, 0) + \sigma(0)A_{-}], \delta\bigg)-\sigma(\delta)A_{+} \Bigg)\\
                    & = x_0 + \delta \mu(x_0, 0) + \frac{\delta}{2}(\sigma(0)A_{-}+\sigma(\delta)A_{+}) + \tilde{\varepsilon}_1,
                \end{split}
            \end{equation}
            where $\tilde{\varepsilon}_1$ satisfies
            \begin{equation}
                \begin{split}
                    \label{eq:leltlitja}
                    \tilde{\varepsilon}_1 &= O(\delta^2)\times(\mu(x_0, 0)+\sigma(0)A_{-})\\
                    &\quad -(\delta/2+O(\delta^3))\\
                    &\quad \times\Bigg([1+O(\delta)]\times\bigg[\mu(x_0, 0)-\mu\Big(x_0+[\delta+O(\delta^2)]\times[\mu(x_0, 0)+\sigma(0)A_{-}], \delta\Big)\bigg]\Bigg)\\
                    &\quad +O(\delta^3)\sigma(\delta)A_{+}.
                \end{split}
            \end{equation}
            Let us define the stochastic process $g$ \yvann{modif}
            \begin{equation}
                g(t) = \mu(x_0+t((\delta+O(\delta^2)))\times(\mu(x_0, 0)+\sigma(0)A_{-}), t\delta),
            \end{equation}
            for $t\in [0, \delta]$.
            Since $\mu$ is uniformly bounded and $A_{-}$ is a normal random variable, $g$ is square-integrable. By~\citep[Theorem 3.2]{CORTES2007757}, there exists a measurable function $\xi\in [0, 1]$ such that
            \begin{equation}
                \begin{split}
                    &\indent \mu(x_0 + (\delta+O(\delta^2))\times(\mu(x_0, 0)+\sigma(0)A_{-}), \delta) - \mu(x_0, 0) \\
                    &= g(1)-g(0)\\
                    &= \dot{g}(\xi)\\
                    &= \delta \partial_t \mu (x_0+\xi (\delta+O(\delta^2))\times(\mu(x_0, 0)+\sigma(0)A_{-}), \xi \delta)\\
                    &+ ((\delta + O(\delta^2))\times (\mu(x_0, 0)+\sigma(0)A_{-}))\partial_x \mu (x_0+\xi (\delta+O(\delta^2))\times(\mu(x_0, 0)+\sigma(0)A_{-}), \xi \delta).
                \end{split}
            \end{equation}
            Let $(\theta, \tau) = (x_0+\xi (\delta+O(\delta^2))\times(\mu(x_0, 0)+\sigma(0)A_{-}), \xi \delta) \in [x_0, x_0 + \delta (\delta+O(\delta^2))\times(\mu(x_0, 0)+\sigma(0)A_{-})]\times [0, \delta]$, then $(\theta, \tau)$ is measurable because $\xi$ is measurable. Furthermore,
            \begin{equation}
                \label{eq:tlek}
                \begin{split}
                    &\quad\mu\bigg(x_0+(\delta+O(\delta^2))\times(\mu(x_0, 0)+\sigma(0)A_{-}),\delta\bigg)-\mu(x_0, 0)\\
                    &=\delta\partial_t \mu(\theta, \tau)+(\delta+O(\delta^2))\times(\mu(x_0,0)+\sigma(0)A_{-})\partial_x\mu(\theta, \tau).
                \end{split}
            \end{equation}
            Since all the involved terms are uniformly bounded, the expectation and $\mathcal{L}^2$-norm of $g(1)-g(0)$ are finite.
            Using $\normLtwo{A_{-}}=O(\delta^{-0.5})$ and $\bbE[A_{-}] = 0$,~\eqref{eq:tlek} is bounded by
            \begin{equation}
                \label{eq:tlek2}
                \begin{split}
                    \bignormLtwo{\mu\big(x_0+(\delta+O(\delta^2))\times(\mu(x_0, 0)+\sigma(0)A_{-}),\delta\big)-\mu(x_0, 0)} &= O(\delta^{0.5}),\\
                    \bbE\big[\mu\big(x_0+(\delta+O(\delta^2))\times(\mu(x_0, 0)+\sigma(0)A_{-}),\delta\big)-\mu(x_0, 0)\big] &= O(\delta).
                \end{split}
            \end{equation}
            Injecting~\eqref{eq:tlek2} into~\eqref{eq:leltlitja}, we have
            \begin{equation}
                \begin{split}
                    \label{eq:noise_eps1}
                    \normLtwo{\tilde{\varepsilon}_1} = O(\delta^{1.5}),\quad \bbE[\tilde{\varepsilon}_1] = O(\delta^2).
                \end{split}
            \end{equation}
            Recall that under Assumption~\ref{assumption:reg_of_sde_flow}, $\mu$ is uniformly bounded and admits uniformly bounded partial derivatives, furthermore, $\sigma$ and $\partial_t\sigma$ are bounded on $[0, T]$ since it is continuously differentiable.
            Thus, there exists $L > 0$ such that uniformly on $\bbR\times [0, T]$
            \begin{equation}
                \label{eq:unif_bounded}
                \max\{\abs{\sigma},\abs{\partial_t\sigma}, \abs{\mu}, \abs{\partial_x\mu},\abs{\partial_t\mu}\}\leq L.
            \end{equation}
            Using~\eqref{eq:unif_bounded}, the $O$-terms derived in~\eqref{eq:tlek2} and~\eqref{eq:noise_eps1} have constants depending on the uniform bound $L$.
            Now, we take into account the sampling error.
            Using~\eqref{eq:kalman_filter_update_eq} and Lemma~\ref{lemma:kalman_gain}, we obtain that the variance of the first component is $P(\delta)_{0,0}=\eta^2\delta^3/12 + O(\delta^5)$, it follows that the sample $\wtX_{\delta}$ satisfies
            \begin{equation*}
                \wtX_{\delta}=m^{(0)}(\delta)+\tilde{\nu}_1,
            \end{equation*}
            with $\tilde{\nu}_1\sim \calN(0, \eta^2\delta^3/12+O(\delta^5))$ and where $m^{(0)}(\delta)$ is the predicted mean~\eqref{eq:local_step_ekf0}.
            This concludes the proof.
        \end{proof}
        \begin{remark}
            Let $R=c\delta^2$ for some constant $c\geq 0$.
            By Remark~\ref{rem:kalman_gain}, the equations~\eqref{eq:local_step_ekf0},~\eqref{eq:leltlitja} and~\eqref{eq:noise_eps1} are still valid, and the above local expansion still holds.
        \end{remark}
    \end{lemma}
    \begin{lemma}[Local step of the parabola-ODE]
        \label{lemma:local_step_parabola_ode}
        Let $\delta > 0$ and let $X$ denote the solution on $[0, \delta]$ of the ODE given by~\eqref{eq:random_ode} with the polynomial approximation given by~\eqref{eq:parabola}.
        Then under Assumption~\ref{assumption:reg_of_sde_flow},
        \begin{equation*}
            X(\delta) = x_0 + \delta\mu(x_0, 0) + \frac{\delta}{2}(\sigma(0)A_{-}+\sigma(\delta)A_{+})+\varepsilon_1,
        \end{equation*}
        where $\varepsilon_1$ satisfies
        \begin{equation*}
            \normLtwo{\varepsilon_1} = O(\delta^{1.5}), \quad \bbE[\varepsilon_1] = O(\delta^2),
        \end{equation*}
        and where the constants in front of the big-$O$ term depend only on uniform bounds on $\mu$ and $\sigma$.
        \begin{proof}
            Applying the mean-value theorem~\citep{CORTES2007757} to $g(t) = \mu(t X(u), tu)$, we obtain that for any $u\in [0, \delta]$, there exist $(\theta_u, \tau_u)$ in $[x_0, X(u)]\times [0, u]$ such that
            \begin{equation}
                \label{eq:mvt_2}
                \mu(X(u),u)= \mu(x_0, 0) + u \partial_t \mu(\theta_u, \tau_u)+(X(u)-x_0)\partial_x \mu(\theta_u, \tau_u).
            \end{equation}
            Injecting~\eqref{eq:mvt_2} into the integrated ODE equation, we obtain
            \begin{equation}
                \label{eq:parabola_ode_one_step}
                \begin{split}
                    X(\delta) &= x_0 + \int\limits_{0}^{\delta}\mu(X(u),u)\dd u + \int\limits_{0}^{\delta}\sigma(u)\dd \beta(u)\\
                    &= x_0 + \delta \mu(x_0, 0) + \int\limits_{0}^{\delta}[u\partial_t \mu(\theta_u, \tau_u)+(X(u)-x_0)\partial_x \mu(\theta_u, \tau_u)]\dd u+\int\limits_{0}^{\delta}\sigma(u)\dd \beta(u).\\
                \end{split}
            \end{equation}
            Using~\eqref{eq:unif_bounded},~\eqref{eq:parabola_ode_one_step}, $\dd \beta (s) = (A_{-}+s\ddot{\beta})\dd s$, and the triangular inequality, we have
            \begin{equation}
                \label{eq:normltwointe}
                \begin{split}
                    \normLtwo{X(u)-x_0} &\leq \biggnormLtwo{\int_{0}^{u}\mu(X(s), s)\dd s +  \int_{0}^{u}\sigma(s)\dd \beta(s))}\\
                    &\leq L\int_{0}^{u}\dd s + L\int_{0}^{u}(\normLtwo{A_{-}}+s \cdot \normLtwo{\ddot{\beta}})\dd s\\
                    &\leq L ( u + \normLtwo{A_{-}} \cdot u+\normLtwo{\ddot{\beta}} \cdot u^2/2).
                \end{split}
            \end{equation}
            Now, let us tackle the second-order term in~\eqref{eq:parabola_ode_one_step}, using~\eqref{eq:unif_bounded} and~\eqref{eq:normltwointe},
            \begin{equation}
                \begin{split}
                    \label{eq:coakokca}
                    &\biggnormLtwo{\int_{0}^{\delta}[u\partial_t \mu(\theta_u, \tau_u)+(X(u)-x_0)\partial_x\mu(\theta_u,\tau_u)]\dd u}\\
                    &\quad \leq L\delta^2/2+L^2(\delta^2/2+\normLtwo{A_{-}} \cdot \delta^2/2+\normLtwo{\ddot{\beta}} \cdot \delta^3/6).
                \end{split}
            \end{equation}
            Let $\varepsilon_1 = X(\delta)-x_0-\delta\mu(x_0)-\frac{\delta}{2}(\sigma(0)A_{-}+\sigma(\delta)A_{+})$.
            Using~\eqref{eq:parabola_ode_one_step} and the triangular inequality, we have
            \begin{equation}
                \label{eq:balbala}
                \begin{split}
                    \normLtwo{\varepsilon_1}&\leq \biggnormLtwo{\int_{0}^{\delta}[u\partial_t \mu(\theta_u, \tau_u)+(X(u)-x_0)\partial_x\mu(\theta_u,\tau_u)]\dd u}
                    \\&\quad+\biggnormLtwo{\int_{0}^{\delta}\sigma(u)\dd \beta(u) - \frac{\delta}{2}(\sigma(0)A_{-}+\sigma(\delta)A_{+})}.
                \end{split}
            \end{equation}
            Using the bounds on $A_{-}$, $\ddot{\beta}$~\eqref{eq:bound_l2_norm}, and~\eqref{eq:coakokca}, the first term of RHS in~\eqref{eq:balbala} is dominated by $O(\delta^{1.5})$ where the constant depends only on uniform bounds on $\mu$ and $\sigma$.
            Using Lemma~\ref{lemma:diff_integral}, the second term of RHS in~\eqref{eq:balbala} is also dominated by $O(\delta^{1.5})$ with a constant depending only on uniform bounds on $\sigma$.
            It follows that
            \begin{equation*}
                \begin{split}
                    \normLtwo{\varepsilon_1}=O(\delta^{1.5}),
                \end{split}
            \end{equation*}
            where the big-$O$ term has a constant depending only on uniform bounds on $\mu$ and $\sigma$, i.e., $L$.
            For the expectation, using $\bbE[\int_{0}^{\delta}\sigma(u)\dd \beta(u)-\frac{\delta}{2}(\sigma(0)A_{-}+\sigma(\delta)A_{+})]=0$, we have
            \begin{equation*}
                \begin{split}
                    \bbE[\varepsilon_1] &= \bbE\bigg[\int_{0}^{\delta}[u\partial_t \mu(\theta_u, \tau_u))+(X(u)-x_0)\partial_x\mu(\theta_u, \tau_u)]\dd u\bigg]\\
                    &= O(\delta^2) + \int_{0}^{\delta}\bbE[(X(u)-x_0)\partial_x\mu(\theta_u,\tau_u)]\dd u\\
                    &= O(\delta^2),
                \end{split}
            \end{equation*}
            where we used $\abs{\partial_t \mu}\leq L$, $\abs{\partial_x \mu}\leq L$, and $\bbE[X(u)-x_0] = O(u)$.
            This concludes the proof.
        \end{proof}
    \end{lemma}

    \subsection{Proof of Theorem~\ref{th:ekf0-scheme2}}\label{app:ekf0-scheme2}
    The proof of Theorem~\ref{th:ekf0-scheme2} follows from, first, deriving the convergence rate of the EKF0 solution $\wtX$ given by Algorithm~\ref{alg:ekf0-scheme2} to the parabola-ODE solution $X(\cdot)$ given by~\eqref{eq:random_ode}, and second, combining this rate to the convergence rate of the parabola-ODE to the true process $X$ given in~\citep[][Th. 3.14., Th. 3.17]{Foster2020OptimalPolyBrownian}.

    The convergence of the EKF0 solution to the parabola-ODE solution is obtained by carefully comparing, in the $\calL^2$-space of square-integrable random variables, the local parabola-ODE expansion and one step of Algorithm~\ref{alg:ekf0-scheme2} with EKF0 linearisation. By doing this, we obtain uniform bounds dominated by $O(\delta^{1.5})$ for the local error and uniform bounds dominated by $O(\delta)$ for the global error. Since the local error rate for the parabola-ODE method is dominated by $O(\delta^2)$ and the global error rate is dominated by $O(\delta)$, we conclude that our method is of order $1.5$ locally and $1.0$ globally.

    \begin{refproof}[Proof of Theorem~\ref{th:ekf0-scheme2}]
        \label{proof:ekf0-scheme2}
        Let $\delta > 0$, $X$ denotes the solution of the piecewise parabola-ODE on $[0, T]$ with $t_n = n\delta$, and $K=T/\delta$. Let $\wtX$ denotes the EKF0 solution given by~\ref{alg:ekf0-scheme2}.
        For any $n\in [1, K]$, we set $E_n = \normLtwo{X(t_n)-\wtX_{t_n}}^2$.
        We aim to prove that $\max_{n\in [1,K]} E_n = O(\delta^2)$.
        For $n=1$, Lemma~\ref{lemma:local_step_parabola_ode} gives
        \begin{equation}
            \label{eq:proof_scheme2_0}
            \begin{split}
                X(\delta) = x_0 + \delta \mu(x_0, 0) + \frac{\delta}{2}(\sigma(0)A_{-,0}+\sigma(\delta)A_{+,0})+\varepsilon_{1},
            \end{split}
        \end{equation}
        where $A_{\pm,0}$ denotes the derivative of $\beta_0$ evaluated at $0$ and $\delta$, and where $\varepsilon_1$ satisfies
        \begin{equation}
            \label{eq:eq_noise}
            \begin{split}
                \normLtwo{\varepsilon_1}^2 = O(\delta^3), \quad \bbE[\varepsilon_1] = O(\delta^2),
            \end{split}
        \end{equation}
        and where the constants in the big-$O$ term depend only on uniform bounds on $\mu$ and $\sigma$.
        Using~\ref{lemma:localekf0}, we have
        \begin{equation}
            \label{eq:proof_global_theorem_1}
            \begin{split}
                \wtX_{\delta} = x_0+\delta\mu(x_0, 0)+\frac{\delta}{2}(\sigma(0)A_{-,0}+\sigma(\delta)A_{+,0})+\tilde{\varepsilon}_1 + \tilde{\nu}_1,
            \end{split}
        \end{equation}
        where $\tilde{\nu}_1\sim\calN(0,\eta^2\delta^3/12+O(\delta^5))$ is independent from $\tilde{\varepsilon}_1$ and $\tilde{\varepsilon}_{1}$ satisfies the same property as $\varepsilon_1$~\eqref{eq:eq_noise}.
        By comparing~\eqref{eq:proof_scheme2_0} and~\eqref{eq:proof_global_theorem_1}, we have
        \begin{equation}
            \label{eq:proof_global_theorem_2}
            \begin{split}
                E_1 &= \bbE[\abs{\varepsilon_1-\tilde{\varepsilon}_1-\tilde{\nu}_1}^2]\\
                &=O(\delta^3).
            \end{split}
        \end{equation}
        This concludes the proof for the first iteration step.
        Let $n\geq 2$.
        Using Lemma~\ref{lemma:local_step_parabola_ode} for the parabola-ODE on $[t_{n-1}, t_n]$, we have the following local expansion
        \begin{equation}
            \label{eq:local_parabola_eq}
            X(t_{n}) = X(t_{n-1})+\delta\mu(X(t_{n-1}), t_{n-1})+\frac{\delta}{2}(\sigma(t_{n-1})A_{-,(n-1)}+\sigma(t_{n})A_{+,(n-1)})+\varepsilon_n,
        \end{equation}
        where $\varepsilon_n$ satisfies~\eqref{eq:eq_noise}.
        Similarly, using Lemma~\ref{lemma:localekf0}, we have
        \begin{equation}
            \label{eq:local_scheme_eq}
            \wtX_{t_{n}} = \wtX_{t_{n-1}}+\delta\mu(\wtX_{t_{n-1}}, t_{n-1})+\frac{\delta}{2}(\sigma(t_{n-1})A_{-,(n-1)}+\sigma(t_{n})A_{+,(n-1)})+\tilde{\varepsilon}_n+\tilde{\nu}_{n},
        \end{equation}
        where $\tilde{\nu}_n$ is the independent sampling error, i.e., $\tilde{\nu}_n\sim\calN(0, \eta^2\delta^3/12+O(\delta^5))$ and $\tilde{\varepsilon}_n$ satisfies~\eqref{eq:eq_noise}.
        The two previous developments lead to
        \begin{equation}
            \label{eq:proof_global_theorem_3}
            \begin{split}
                E_n&=\bbE[\abs{X(t_{n-1})-\wtX_{t_{n-1}}\\
                &\quad +\delta(\mu(X(t_{n-1}),t_{n-1})-\mu(\wtX_{t_{n-1}},t_{n-1}))\\
                &\quad +\varepsilon_n-\tilde{\varepsilon}_n-\tilde{\nu}_n}^2].\\
            \end{split}
        \end{equation}
        Since $\mu$ has bounded partial derivatives, $\mu$ is uniformly Lipschitz with constant $L$~\eqref{eq:unif_bounded}.
        Expanding~\eqref{eq:proof_global_theorem_3}, using the Lipschitzness of $\mu$ and the independence of $\nu_n$ from the other variables, we have
        \begin{equation}
            \label{eq:proof_global_theorem4}
            \begin{split}
                E_n &\leq (1 + 2L\delta + L^2\delta^2)E_{n-1}\\
                &\quad + 2\bbE[(\varepsilon_n-\tilde{\varepsilon}_n)(X(t_{n-1})-\wtX_{t_{n-1}}+\delta(\mu(X(t_{n-1}),t_{n-1})-\mu(\wtX_{t_{n-1}},t_{n-1})))] \\&\quad + \bbE[|\varepsilon_n-\tilde{\varepsilon}_n|^2]+\bbE[|\tilde{\nu}_n|^2].
            \end{split}
        \end{equation}
        Similarly to~\eqref{eq:proof_global_theorem_2}, the last terms are dominated by $O(\delta^3)$,
        \begin{equation*}
            \bbE[|\varepsilon_n-\tilde{\varepsilon}_n|^2]+\bbE[|\tilde{\nu}_n|^2] = O(\delta^3).
        \end{equation*}
        We need to control the expectation term $\bbE[(\varepsilon_n-\tilde{\varepsilon}_n)(X(t_{n-1})-\wtX_{t_{n-1}}+\delta(\mu(X(t_{n-1}),t_{n-1})-\mu(\wtX_{t_{n-1}},t_{n-1})))]$.
        Using~\eqref{eq:local_parabola_eq} and~\eqref{eq:local_scheme_eq} for $r\leq n-1$, we have
        \begin{equation*}
            \begin{split}
                X(t_{n-1})-\wtX_{t_{n-1}} &= \sum_{r=1}^{n-1}\varepsilon_r-\tilde{\varepsilon}_r-\tilde{\nu}_r+\delta(\mu(X(t_{r-1}), t_{r-1})-\mu(\wtX_{t_{r-1}}, t_{r-1})).
            \end{split}
        \end{equation*}
        Hence, we need to control the following terms
        \begin{equation}
            \label{eq:eq_proof_global_theorem7}
            \begin{split}
                \bbE\bigg[(\varepsilon_n-\tilde{\varepsilon}_n)\bigg(\sum_{r=1}^{n-1}\varepsilon_r-\tilde{\varepsilon}_r-\tilde{\nu}_r\bigg)\bigg],\quad \bbE\bigg[(\varepsilon_n-\tilde{\varepsilon}_n) \bigg(\sum_{r=1}^{n-1}\mu(X(t_{r}), t_{r})-\mu(\wtX_{t_{r}}, t_{r})\bigg)\bigg].
            \end{split}
        \end{equation}
        Let $j\in [1, K]$.
        Using~\eqref{eq:parabola_ode_one_step}, we obtain for some $(\theta_u^{(j)},\tau_u^{(j)})\in [X(t_{j-1}), X(t_j)]\times [t_{j-1}, t_j]$,
        \begin{equation*}
            \begin{split}
                \varepsilon_j &= \underbrace{\int_{t_{j-1}}^{t_j}[u\partial_t\mu(\theta_u^{(j)},\tau_u^{(j)})+(X(u)-X(t_{j-1}))\partial_x\mu(\theta_u^{(j)}, \tau_u^{(j)})]\dd u}_{R_{1,j}} \\&\quad+ \underbrace{\int_{t_{j-1}}^{t_j}\sigma(u)\dd \beta(u) - \frac{\delta}{2}(\sigma(t_{j-1})A_{-,j-1}+\sigma(t_j)A_{+,j-1})}_{R_{2,j}},
            \end{split}
        \end{equation*}
        with $R_{1,j}$ which can be further expanded,
        \begin{equation}
            \label{eq:expression_R1j}
            \begin{split}
                R_{1,j} &= \int\limits_{t_{j-1}}^{t_j}\bigg[u\partial_t \mu(\theta_u^{(j)}, \tau_u^{(j)})+\bigg(\int\limits_{t_{j-1}}^{u}\mu(X(s), s)\dd s+\int\limits_{t_{j-1}}^{u}\sigma(s)\dd \beta(s)\bigg)\partial_x \mu(\theta_u^{(j)}, \tau_u^{(j)})\bigg]\dd u
                \\ & =O(\delta^2) + \int\limits_{t_{j-1}}^{t_j}\bigg(\int\limits_{t_{j-1}}^{u}\sigma(s)\dd \beta(s)\bigg)\partial_x\mu(\theta^{(j)}_u, \tau^{(j)}_u)\dd u.
            \end{split}
        \end{equation}
        Using~\eqref{eq:expression_R1j}, the independence of the Brownian parabolas between $[t_{r-1}, t_{r}]$ and $[t_{n-1}, t_{n}]$ for $r<n$, and $\bbE[\dd \beta] = 0$, we have $\bbE[R_{1,n}R_{1,r}]=O(\delta^4)$.
        Using the independence of the Brownian parabolas and $\bbE[R_{2,j}]=0$, we have $\bbE[R_{2,n}R_{2,r}]=\bbE[R_{2,n}]\bbE[R_{2,r}]=0$.
        Similarly, we have $\bbE[R_{1,n}R_{2, r}]=0$ and $\bbE[R_{2,n}R_{1,r}]=0$.
        Hence, $\bbE[\varepsilon_n\varepsilon_r] = O(\delta^4)$.
        Similarly, using~\eqref{eq:leltlitja} and~\eqref{eq:tlek}, the independence of two distinct Brownian parabolas and the uniform bounds on $\mu$ and $\sigma$, we have $\bbE[\tilde{\varepsilon}_n\tilde{\varepsilon}_r] = O(\delta^4)$.
        Likewise, it can be shown that $\bbE[\varepsilon_n\tilde{\varepsilon}_r]=O(\delta^4)$, $\bbE[\tilde{\varepsilon}_n\varepsilon_r] = O(\delta^4)$, and $\bbE[\varepsilon_n \tilde{\nu}_r] = \bbE[\tilde{\varepsilon}_n\tilde{\nu}_r] = 0$.
        Hence, the first term in~\eqref{eq:eq_proof_global_theorem7} is dominated by $O(n\delta^4)$,
        \begin{equation}
            \label{eq:premier_terme_difficile}
            \bbE\bigg[(\varepsilon_n-\tilde{\varepsilon}_n)\big(\sum_{r=1}^{n-1}\varepsilon_r-\tilde{\varepsilon}_r-\tilde{\nu}_r\big)\bigg] = O(\delta^4n).
        \end{equation}
        Now, let us tackle the second term in~\eqref{eq:eq_proof_global_theorem7}.
        Using Cauchy-Schwarz and the Lipschitzness of $\mu$, we have
        \begin{equation}
            \begin{split}
                \label{eq:second_terme_difficile}
                \bbE\bigg[\varepsilon_n\bigg(\sum_{r=1}^{n-1}\mu(X(t_{r}), t_{r})-\mu(\wtX_{t_{r}}, t_{r})\bigg)\bigg]&\leq \sqrt{\bbE[\abs{\varepsilon_n}^2]}\sum_{r=1}^{n-1}\sqrt{\bbE[\abs{\mu(X(t_{r}), t_{r})-\mu(\wtX_{t_{r}}, t_{r})}^2]}\\&
                \leq L\sqrt{\bbE[\abs{\varepsilon_n}^2]}\sum_{r=1}^{n-1}\sqrt{\bbE[\abs{X(t_{r})-\wtX_{t_{r}}}^2]}\\
                &= O(\delta^{1.5})\sum_{r=1}^{n-1}\sqrt{E_{r}}.
            \end{split}
        \end{equation}
        Injecting~\eqref{eq:premier_terme_difficile} and~\eqref{eq:second_terme_difficile} into~\eqref{eq:proof_global_theorem4}, we have for a well-chosen constant in front of all the previous $O$-terms, $A$, (which is possible since all the previous $O$-terms have constants depending only on uniform bounds on $\mu$ and $\sigma$)
        \begin{align}
            \label{eq:recursive_error}
            \begin{split}
                E_n&\leq (1+2L\delta+L^2\delta^2)E_{n-1} + 2A\delta^4n + 2A\delta^{2.5}\sum_{r=1}^{n-1}\sqrt{E_{r}} + A\delta^3\\
                &\leq (1+B\delta)E_{n-1} + C\delta^3 + C\delta^{2.5}\sum_{r=1}^{n-1}\sqrt{E_{r}},
            \end{split}
        \end{align}
        with $B=2L+L^2$, $C=3A$ and where we used $\delta \leq 1$, $n\leq T/\delta$.
        By recursively applying~\eqref{eq:recursive_error} with initialisation $E_0 = 0$, we obtain
        \begin{align}
            \label{eq:last_bound_global_theorem}
            E_n = O(n\delta^3).
        \end{align}
        Setting $n=K=T/\delta$ in~\eqref{eq:last_bound_global_theorem} leads to $E_{K} = O(\delta^2)$. Hence, the EKF0 solution $\wtX$ converges to the parabola-ODE solution $X(\cdot)$, locally with order $1.5$ and globally with order $1.0$.
        Using the triangular inequality, we combine the previous convergence rates with the convergence rates of the parabola-ODE to the true process~\citep[][Th. 3.14., Th. 3.17]{Foster2020OptimalPolyBrownian} which are, locally of order $2.0$ and globally of order $1.0$. This concludes the proof.
        \begin{remark}
        \label{rem:convth_with_non_zero_error_measurement}
            By Remark~\ref{rem:kalman_gain}, the previous proof is still valid when $R=c\delta^2$ for some $c\geq 0$.
            Thus, the EKF0 solution given by Algorithm~\ref{alg:ekf0-scheme2} with measurement error $R=O(\delta^2)$ converges with rates as stated above.
        \end{remark}
    \end{refproof}

    \subsection{Proof of Theorem~\ref{th:ekf0-scheme3}}\label{app:ekf0-scheme3}

    The proof of Theorem~\ref{th:ekf0-scheme3} is similar to the proof of Theorem~\ref{th:ekf0-scheme2} but takes into account non-zero variances $(P_k(t_k))$ for the prior distribution.

    \begin{refproof}[Proof of Theorem~\ref{th:ekf0-scheme3}]
        \label{proof:ekf0-scheme3}
        We adapt the previous proof to take into account non-zero variance initialisations $(P_n(t_n))$ and no sampling error throughout the EKF passes.
        Let $\delta > 0$ and $n\geq 1$.
        Let $Y_{n-1}(t_{n})\sim\calN(m_{n-1}(t_{n}), P_{n-1}(t_{n}))$ be the posterior distribution given by Algorithm~\ref{alg:ekf0-scheme3} for the $n$-th EKF0 pass.
        Let $\wtX_{t_{n}}\sim \calL(Y_{n-1}^{(0)}(t_{n}))$.
        Let $E_{n, (1)} = \normLtwo{X(t_{n})-m_{n-1}^{(0)}(t_{n})}^2$ and $E_{n} = \normLtwo{X(t_{n})-\wtX_{t_{n}}}$.
        Using the predictive equations~\eqref{eq:pred_eq_kalman} and~\eqref{eq:dev_Q}, we obtain that for any $n\in [1, K]$
        \begin{equation}
            \label{eq:pred_var_2}
            P^{-}_{n-1}(t_{n}) = \eta^2
            \begin{pmatrix}
                n(\delta^3/3 + O(\delta^4)) & \delta^2/2 + O(\delta^2) \\
                \delta^2/2 + O(\delta^2)    & \delta + O(\delta^2)
            \end{pmatrix},
        \end{equation}
        where $P^{-}_{n-1}$ is the predictive variance of the $n$-th EKF0 pass.
        Using the update equations~\eqref{eq:kalman_filter_update_eq} and adapting Lemma~\ref{lemma:kalman_gain} for non-zero variance $P_{n}(t_n)$, we obtain that for any $n\in [0, K]$
        \begin{equation}
            \label{eq:gain_k2}
            K_n(\delta) = \begin{pmatrix}
                              \delta/2 + O(\delta^3) \\
                              1 + O(\delta^3)
            \end{pmatrix}.
        \end{equation}
        Hence, the local update equation for the mean~\eqref{eq:local_step_ekf0} is still valid.
        Following the proof in~\eqref{proof:ekf0-scheme2} but without adding a sampling error, we obtain for $n\geq 1$
        \begin{equation}
            \label{eq:proof_expr2_scheme3}
            E_{n,(1)} = O(n\delta^3),
        \end{equation}
        where $A$ is a well-chosen uniform constant as in~\eqref{eq:last_bound_global_theorem}.
        Using~\eqref{eq:pred_var_2},~\eqref{eq:gain_k2}, and~\eqref{eq:kalman_filter_update_eq}, we obtain $P_{n-1}(t_{n})\sim \calN(0, n\eta^2(\delta^3/12+O(\delta^5))$ for $n\geq 1$.
        Using~\eqref{eq:proof_expr2_scheme3}, we have for $n\geq 1$
        \begin{align}
            \label{eq:last_line_proof_scheme3}
            E_{n} = \normLtwo{X(n\delta)-\wtX_{t_{n}}}^2 \leq E_{n,(1)}+ n\eta^2(\delta^3/12+O(\delta^5)) = O(n \delta^3).
        \end{align}
        Setting $n=K=T/\delta$ in~\eqref{eq:last_line_proof_scheme3} leads to the same previous error rate, $E_K = O(\delta^2)$. We conclude the proof just as previously.
        \begin{remark}
            Similarly to Remark~\ref{rem:convth_with_non_zero_error_measurement} for Theorem~\ref{th:ekf0-scheme2}, the previous convergence rates still holds when $R=O(\delta^2)$.
        \end{remark}
    \end{refproof}

    \subsection{Marginalised Gaussian SDE Filter in the multivariate case}
    In this section, we describe the state-space model in the multivariate case with the parabola's coefficients as part of the state.
    \label{appendix:multivariate_case_marginal_brownian}
    Let
    \begin{equation}
        Y_{k}^{(i)} = (Y_k^{(i), (0)}, Y_k^{(i), (1)}, W^{(i)}_{k\delta, (k+1)\delta}, I^{(i)}_{k\delta, (k+1)\delta})^{\top}
    \end{equation}
    and let $Y_k = (Y_k^{(1)}, \ldots Y_k^{(d)})^{\top}$ be the concatenation of the univariate states.
    The initial distribution is given by
    \begin{equation*}
        Y_k(t_k)\sim \calN(m_k(t_k), P_k(t_k)),
    \end{equation*}
    where $Y_k(t_k)$ is initialised with
    \begin{equation*}
        \begin{split}
            Y_k^{(i), (0)}(t_k) & \coloneqq\wtX_{t_k}^{(i)} \sim \calL(Y_{k-1}^{(i), (0)}(t_k)),\\
            Y_k^{(i), (1)}(t_k) & \coloneqq f_k^{(i)}(\wtX_{t_k}, t_k)\\
            &=\mu^{(i)}(\wtX_{t_k}, t_k)+\sum_{j=1}^{d}\sigma(t_k)_{i,j}\left(\frac{W^{(j)}_{k\delta,(k+1)\delta}}{\delta}-\frac{\sqrt{6}I^{(j)}_{k\delta,(k+1)\delta}}{\delta}\right),
        \end{split}
    \end{equation*}
    where $B_{k\delta,(k+1)\delta}$ and $I_{k\delta,(k+1)\delta}$ are central normal variables with variance, $\delta I_{d}$ and $\delta/2 I_{d}$, respectively.
    Hence, $m_k(t_k) = (\wtX_{t_k}^{(i)}, \mu^{(i)}(\wtX_{t_k}, t_k), 0, 0)^{\top}$ and $P_k(t_k) = [V_k^{(i,j)}]$, where $V_k^{(i,j)}$ are given by, for $i,j\in [1, d]$
    \begin{equation*}
        \begin{split}
            V_k^{(i,j)} &= \bbC[Y_k^{(i)}(t_k), Y_k^{(j)}(t_k)]\\
            &= \begin{pmatrix}
                   0 & 0                                            & 0                      & 0                                \\
                   0 & \frac{4}{\delta} (\sigma\sigma^{\top})_{i,j} & \sigma_{i,j}           & -\frac{\sqrt{6}}{2}\sigma_{i,j}  \\
                   0 & \sigma_{j,i}                                 & \mathds{1}_{i=j}\delta & 0                                \\
                   0 & -\frac{\sqrt{6}}{2}\sigma_{j,i}              & 0                      & \mathds{1}_{i=j}\frac{\delta}{2}
            \end{pmatrix}(t_k).
        \end{split}
    \end{equation*}
    Assuming the same prior on each coordinate, the transition and variance matrices, $\bar{A}$ and $\bar{Q}$ are given by
    \begin{equation}
        \bar{A} = I_d \otimes \textup{Diag}(A, 1, 1),\quad \bar{Q} = I_d\otimes \textup{Diag}(Q, 0, 0),
    \end{equation}
    where $A$ and $Q$ are the matrices arising in the Markov transition for one coordinate.
    Let $H_0 = I_{d} \otimes \begin{pmatrix}
                                 1 & 0 & 0 & 0
    \end{pmatrix}$, $H_1 = I_{d} \otimes \begin{pmatrix}
                                             0 & 1 & 0 & 0
    \end{pmatrix}$ and $H_{2, 3} = I_{d} \otimes \begin{pmatrix}
                                                     0 & 0 & 1 & 1
    \end{pmatrix}$.
    Let us denote by $H_k$
    \begin{equation*}
        H_k(t)= I_{d} \otimes \begin{pmatrix}
                                  0 & 0 & 0        & 0                                                                \\
                                  0 & 0 & 0        & 0                                                                \\
                                  0 & 0 & 1/\delta & 0                                                                \\
                                  0 & 0 & 0        & \frac{\sqrt{6}}{\delta} \left(\frac{2(t-t_k)}{\delta} - 1\right)
        \end{pmatrix},
    \end{equation*}
    such that $(H_{2,3} H_k(t) Y_k(t))^{(j)} = \frac{B^{(j)}_{k\delta,(k+1)\delta}}{\delta}+\frac{\sqrt{6}I^{(j)}_{k\delta,(k+1)\delta}}{\delta}\left(\frac{2(t-t_k)}{\delta}-1\right)$ for $j\in [1, d]$.
    The vector field $f_k$ satisfies $f_k(Y_{k}^{(0)}(t), t) = \bar{f}_k(Y_{k}(t), t)$ where
    \begin{equation*}
        \bar{f}_k(Y, t) = \mu(H_0 Y, t) + \sigma(t) H_{2,3}H_{k}(t)Y.
    \end{equation*}
    Thus, the EKF0 scheme is given by
    \begin{equation*}
        \tilde{H}_{k}(t) = H_1 - \sigma(t) H_{2,3}H_{k}(t),
    \end{equation*}
    while for EKF1, we replace $\mu(H_0 Y, t)$ by its linearisation around the predicted mean, which leads to
    \begin{equation*}
        \tilde{H}_k(t) = H_1 - (J_{\mu(\cdot, t)}(H_{0}m^{-}(t), s)H_{0} + \sigma(t) H_{2, 3}H_{k}(t)).
    \end{equation*}

    \bibliographystyle{ba}
    \bibliography{main}
\end{document}